\documentclass[11pt,english]{article}
\usepackage[sc]{mathpazo}
\usepackage{geometry}
\usepackage{color}
\usepackage{array}
\usepackage{bbding}
\usepackage{multirow}
\usepackage[fleqn]{amsmath}
\usepackage{amssymb}
\usepackage{amsthm}
\usepackage{graphicx}
\usepackage{natbib}
\usepackage{lscape}
\usepackage[hyphens]{url}
\usepackage[hidelinks]{hyperref}
\usepackage{cleveref}
\usepackage{enumitem}
\usepackage{comment}
\usepackage{bm}
\usepackage[title]{appendix}
\usepackage{longtable}
\usepackage{mathtools}
\usepackage{chngcntr}
\usepackage{authblk}
\usepackage{babel}
\usepackage{dsfont}
\usepackage{subcaption}
\usepackage[most]{tcolorbox}
\usepackage{arydshln}
\usepackage[title]{appendix}
\usepackage{algorithm}
\usepackage{algorithmic}

\usepackage{color-edits}
\addauthor[VS]{vs}{magenta}

\makeatletter
\allowdisplaybreaks
\usepackage[mathlines, pagewise]{lineno}
\usepackage{etoolbox} 

\newcommand*\linenomathpatchAMS[1]{%
	\expandafter\pretocmd\csname #1\endcsname {\linenomathAMS}{}{}%
	\expandafter\pretocmd\csname #1*\endcsname{\linenomathAMS}{}{}%
	\expandafter\apptocmd\csname end#1\endcsname {\endlinenomath}{}{}%
	\expandafter\apptocmd\csname end#1*\endcsname{\endlinenomath}{}{}%
}

\expandafter\ifx\linenomath\linenomathWithnumbers
\let\linenomathAMS\linenomathWithnumbers
\patchcmd\linenomathAMS{\advance\postdisplaypenalty\linenopenalty}{}{}{}
\else
\let\linenomathAMS\linenomathNonumbers
\fi

\linenomathpatchAMS{gather}
\linenomathpatchAMS{multline}
\linenomathpatchAMS{align}
\linenomathpatchAMS{alignat}
\linenomathpatchAMS{flalign}

\newtheorem{prop}{Proposition}
\newtheorem{cor}{Corollary}[prop]

\DeclareMathOperator*{\argmax}{arg\,max}
\begin{document}

\title{Managing ride-sourcing drivers at transportation terminals: a lottery-based queueing approach}
\author[1]{Tianming Liu}
\author[1]{Yafeng Yin}
\author[2]{Vijay Subramanian}

\affil[1]{\small\emph{Department of Civil and Environmental Engineering, University of Michigan, Ann Arbor, MI, United States}\normalsize}
\affil[2]{\small\emph{ 
Department of Electrical Engineering and Computer Science, University of Michigan, Ann Arbor, MI, United States}\normalsize}

\date{September 15, 2025}
\maketitle

\begin{abstract}
\noindent \textbf{Problem definition:} Transportation terminals such as airports often experience persistent oversupply of idle ride-sourcing drivers, resulting in long driver waiting times and inducing externalities such as curbside congestion. While platforms now employ virtual queues with control levers like dynamic pricing, information provision, and direct admission control to manage this issue, all existing levers involve significant trade-offs and side effects. This limitation highlights the need for an alternative management approach. \textbf{Methodology/results:} We develop a queueing-theoretic framework to model ride-sourcing operations at terminals and propose a novel lottery-based control mechanism for the virtual queue. This non-monetary strategy works by probabilistically assigning a driver's entry position. By directly influencing their expected waiting time, the mechanism in turn shapes their decision to join the queue. We reformulate the resulting infinite-dimensional, non-smooth optimization into a tractable bi-level program by leveraging the threshold structure of the equilibrium. Theoretically, we prove that the lottery mechanism can achieve higher or equal social welfare than FIFO-queue-based dynamic pricing. Numerical experiments in unconstrained markets show that in profit maximization, our approach only narrowly trails dynamic pricing and significantly outperforms static pricing. Furthermore, it is shown that under commission fee caps, the lottery mechanism can surpass dynamic pricing in profitability. \textbf{Implications:} This study introduces a new, non-monetary lever for managing idle ride-sourcing drivers at transportation terminals. By aligning operational practices with queue-based dynamics, the proposed lottery mechanism offers a robust and implementable alternative to pricing-based approaches, with advantages in both unconstrained and regulated markets.
\end{abstract}

\section{Introduction}\label{sec:Intro}

Transportation terminals, such as airports and railway stations, often face a large number of arriving travelers who need local last-mile service to go to their destinations. Nowadays, ride-sourcing services, such as Uber, Lyft, and DiDi Chuxing, have captured a significant portion of the terminal passenger demand. According to reports \citep{SFOreport,BALreport}, ride-sourcing has become the largest passenger service at the San Francisco International Airport and Ronald Reagan Washington National Airport, occupying 33\% and 37\% of the market share among all ground transportation for passengers. Consequently, the large concentration of idle ride-sourcing drivers at transportation terminals creates challenges for both platforms and operators. For the platforms, this creates localized oversupply \citep{xu2020empty,pollio2021uber,liu2025bounded}, exacerbating the spatial imbalance of supply and demand, and leading some drivers to wait for hours for a single trip \citep{NYCHearing}. Simultaneously, terminal operators must contend with the externalities of this vehicle influx, which include severe curbside and roadway congestion \citep{hermawan2018impacts,leiner2020transportation}, diminished parking availability for the public \citep{wadud2020examination}, and increased local emissions \citep{agarwal2023information}.

To mitigate the impact of idle drivers on the general traffic in terminal areas, ride-sourcing platforms now apply virtual queues \citep{pollio2021uber,UberHKG,Lyftairport} to manage their operations at these locations. Generally, the platforms establish an exclusive staging area where idle drivers can join a virtual queue and wait to be matched to a customer by the platforms. Therefore, for efficient terminal ride-sourcing operations, management of these ride-sourcing virtual queues has itself become a critical operational challenge for platforms, spurring strategies from both academic literature and industry practice. The first approach is dynamic pricing, in which the platform can dynamically adjust the commission rate of the ride based on the condition of the queue. This flexibility grants the platform a direct mechanism for balancing driver supply with passenger demand. For a queue operating under a first-in-first-out (FIFO) discipline, this can achieve profit and social-welfare optimality \citep{chen2001state}. A second approach is information provision, in which the platform broadcasts information to drivers to steer their behavior toward desired outcomes. For example, \cite{ji2021automated} showed that providing selected drivers who are on their way to terminals with terminal demand predictions can proactively reduce oversupply. Similarly, \cite{agarwal2023information} found that offering real-time queue data allowed drivers to self-distribute more effectively between different terminal locations, thereby reducing supply imbalances. Apart from these strategies, platforms also employ methods of direct admission control. This strategy, commonly observed in current industry practice, involves setting a firm capacity limit for the queue. Once the queue is full, all newly arriving drivers are rejected from entering until a space opens up \citep{pollio2021uber,adarkwah2023reducing}. By directly controlling the number of drivers in the virtual queue, this method provides a straightforward safeguard against excessive oversupply and the resulting congestion. Recent studies have proposed joint strategies that combine multiple approaches for virtual queue management. For example, \cite{varma2023dynamic} proposed applying dynamic pricing-matching and providing long-term average waiting time to drivers together, and proved their proposed scheme can achieve near-optimal control in a probabilistic fluid model of virtual queues. Building on this, \cite{yang2024learning} utilized learning-based dynamic pricing with priority matching to prove near-optimality even under queues with unknown passenger arrival rates.

While the established terminal queue management approaches have all been proven effective in certain respects, each of them also has well-documented trade-offs. Dynamic pricing, while effective at improving system performance, can be unpopular with drivers and has raised broader concerns about fairness and transparency \citep{ashkrof2020understanding, seele2021mapping}. Information provision may trigger undesirable side effects that worsen oversupply. It has been empirically demonstrated that, controlling for existing queue conditions, providing wait time or demand information can increase a driver's propensity to join the queue \citep{liu2025bounded} and prolong their waiting tolerance \citep{xu2020empty,yu2022delay}. Direct admissions control, meanwhile, has been observed not to be able to effectively turn drivers away due to the lack of incentives. Empirical studies have documented that when rejected from a full virtual queue, drivers simply relocate to nearby roads and curbsides waiting for an opening. This behavior, in turn, creates significant traffic interference and external congestion \citep{adarkwah2023reducing,SFOqueue}. In summary, these limitations and side effects of existing approaches motivate the need for alternative strategies.

To this end, we propose and analyze a lottery-based mechanism as an alternative queue management approach for terminal ride-sourcing. This lottery-based mechanism manages driver waiting times by probabilistically assigning them to different starting positions within the virtual queue, offering an alternative to existing strategies. To analyze this mechanism, we model this system using a
queueing framework with state-dependent arrival rates and unlimited capacity, where heterogeneous drivers
make rational decisions on whether to enter the terminal virtual queue. To address the inherent complexities of the queueing dynamics in our model, we characterize the structural properties of the optimal equilibrium to enable scalable optimization of the lottery control policy. Theoretically, we prove that the proposed lottery control's superiority against dynamic pricing in social welfare maximization. Our numerical analysis further demonstrates the robust and competitive performance of the proposed approach under various market conditions, regulations, and objectives.

Our paper makes the following contributions:
\begin{itemize}

    \item We formulate a novel, non-monetary control mechanism for ride-sourcing queues at terminals. This lottery-based approach uses probabilistic assignments of queue positions to influence driver participation and enhance system-wide efficiency.

    \item We theoretically prove our lottery control is superior to dynamic pricing for maximizing social welfare. In profit maximization, numerical analysis confirms our scheme's robust performance, showing it nearly matches dynamic pricing without commission caps and surpasses it when such caps are enforced, while significantly outperforming static pricing.

    \item We develop a solution framework for the optimal lottery control problem, which is naturally an infinite-dimensional and non-smooth optimization. By analyzing the structural properties of the market equilibrium, we reformulate this intractable problem into a well-behaved bi-level program. We prove that this transformation preserves the social welfare optimum, enabling efficient and scalable implementation.

\end{itemize}

The remainder of this paper is organized as follows: \Cref{sec:LR} reviews the relevant literature. \Cref{sec:model} introduces the queueing system model of terminal ride-sourcing operations, explains the lottery control, and formulates the lottery control optimization problems. \Cref{sec:theory} presents our theoretical analysis of the equilibrium properties of the optimal solutions and our proposed problem transformation as well as the transformation's properties. Numerical studies on the solution properties and sensitivity analysis showcase the lottery control's properties and its relative performance compared to pricing strategies are illustrated in \Cref{sec:numerical}. At the end, \Cref{sec:conclusion} concludes the paper.

\section{Related literature} \label{sec:LR}
Aiming to manage the behavior of idle drivers to improve ride-sourcing system performance through a lottery in a virtual queue setting, our work is highly relevant to the literature on decentralized control of strategic drivers in ride-sourcing systems. In this section, we review the established literature, organized by control levers and modeling approaches, and position our work.

A large amount of work focuses on dynamic fare pricing, such as the well-known surge pricing scheme \citep{zha2018surge}, to conduct decentralized driver control. By deterring/attracting supply and demand through the instrument of trip fare pricing,  the platform can improve the system performance by reducing the imbalance between demand and supply. Analytical research on this topic has primarily followed two modeling paradigms. The first uses aggregated queueing models to capture system-wide dynamics. For example, \cite{freund2021pricing} employed a multi-class queueing model to show that dynamic pricing is crucial for stabilizing the ride-sourcing market and preventing oscillations between oversupply and undersupply. The second approach uses leader-follower games and network flow models to analyze strategic interactions. Using this approach, a series of works \citep{bimpikis2019spatial,guda2019your,besbes2021surge,hu2022surge} have established the effectiveness of dynamic fare pricing on service balancing and system efficiency improvement. Other studies in this area have used models such as non-equilibrium dynamic systems \citep{nourinejad2020ride} and simulation \citep{chen2019inbede,chen2021spatial} to explore numerical solution methods of optimal fare pricing.

The aforementioned studies generally assume that the matching between drivers and customers is random or FIFO. To more strategically match drivers with customers, some studies pair customized matching with dynamic fare pricing to improve the system performance. On the local level, \cite{varma2023dynamic} studied the joint optimal pricing and matching between heterogeneous customers and drivers using an open queue model. Assuming that the long-term average waiting time is broadcast to the drivers, they proved that a joint dynamic pricing and max weight matching scheme can achieve near-optimal performance in profit maximization. Assuming that the demand and supply's arrival rates are unknown, \cite{yang2024learning} proposed an online-learning-based dynamic fare pricing model with a static priority matching scheme for an open, multi-class queueing model, and proved that their scheme can achieve near-optimum profit in the long run. On the network level, \cite{ozkan2020joint} revealed that the joint adaptation of dynamic pricing and matching is critical for robust optimization of system performance using a network flow model. \cite{xu2021generalized} introduced a generalized fluid model featuring a flexible, general matching function. They used this framework to develop algorithms for jointly optimizing pricing and matching across the entire ride-sourcing network.


A second lever involves other monetary incentives that offer drivers a direct payment for undertaking a desired action without altering the price on the passenger side. Such direct monetary incentives alter a driver's movement decision by offering a targeted subsidy that artificially increases the profitability of moving to a desired area. These incentives can take two main formats: relocation subsidies, which provide an upfront bonus to drivers for moving to a specific area \citep{zhu2021mean,wangexpress}, and guaranteed wages, which offer a safety net by providing a minimum payment if drivers fail to secure a ride after relocating \citep{xu2020recommender,yengejeh2021rebalancing,zhang2024dynamic}. Theoretical research has established the existence of an optimal subsidy set using network flow models under a leader-follower game setting \citep{wangexpress}. Computationally, various methods have been explored to solve for optimal policies in complex ride-sourcing models. These approaches include simulation-based optimization \citep{yengejeh2021rebalancing}, approximation algorithms \citep{zhang2024dynamic}, and reinforcement learning \citep{zhu2021mean}, which are applied to simulation or mean-field game models.

Apart from monetary incentives, non-monetary levers have also seen applications in the decentralized control problem. In this approach, the platform provides information about expected demand/demand-supply-balance/earnings in selected locations to drivers, hoping to influence their strategic reposition behavior by curating the data presented to drivers, thereby shaping their expectations of future opportunities in different locations. On this approach, studies have mostly utilized network flow models to capture the movement of drivers in the system, and explored strategies including providing drivers with system states like global supply distribution \citep{sadeghi2019re} or customized heatmaps \citep{alnaggar2024heatmap, haferkamp2024heatmap}. These studies also devised optimization approaches to determine the ideal information design. In another study, \cite{beojone2023relocation} modeled the drivers' strategic movements using a Markov decision process, incorporating future revenue forecasts to optimize their long-term earnings.

Our paper differs from the existing literature in two primary ways. First, existing studies, regardless of their modeling scope, generally utilize a long-term aggregate waiting time independent of queue status based on random matching, FIFO matching, or spatial matching to determine driver supply. However, as terminals are high-demand and attractive areas for ride-sourcing drivers, the queue condition will fluctuate significantly across time,  which heavily influences a driver's decision to join or leave and consequently affects the platform operations. We account for this issue by tracking the state-dependent waiting time and the drivers' responses to it in our model. Second, we utilize a different kind of control lever from existing studies. Instead of monetary levers or information levers, we utilize a probabilistic queue entry position lottery to impact the behavior of strategic drivers through controlling waiting time.


\section{Model and Problem Formulation} \label{sec:model}

\subsection{Setting}
We consider a ride-sourcing system, which contains a transportation terminal area and an outside area. On the demand side, we assume that ride-sourcing passengers arrive at the terminal with a Poisson process of rate $\mu$. On the supply side, we assume that the drivers in the system can be classified into $M$ heterogeneous groups based on their valuation of completing an order from the terminal area. Group $m$ ($m\in\{1,2,\cdots,M\}$) drivers arrive at the transportation terminal each following an independent Poisson process of rate $\lambda_m$.

The ride-sourcing platform utilizes a virtual queue for idle drivers at a terminal. To receive passenger matches, drivers must be in the virtual queue and stay within the terminal area. Any arriving passenger is matched with the driver at the front of the queue. An idle group-$m$ driver arriving at the terminal decides whether to join based on the current queue length, $n$. The platform provides the expected waiting time, $W_m^n$, and a commission fee, $p_m$. If the driver joins, they wait to be matched. After a match, the driver completes the trip, which takes an average duration of $T_d$. Thus, joining the queue means that the driver would have to spend an expected time of $W_m^n+T_d$ in the queue and delivering the passenger. The expected utility $U_{m}^n$ that the driver can obtain from joining the queue at length $n$ is:
\begin{equation}
    U_{m}^n=R_m-p_m-r(T_d+W_m^n)
    \label{eq:join_utility}
\end{equation}
\noindent in which $R_m$ is the utility obtained by a group-$m$ driver from the terminal order; and $r$ is the utility-gaining rate in the general market outside of the terminal area. $R_m-p_m$ is the net revenue of the driver gains from joining the queue, while $r(T_d+W_m^n)$ is the incurred opportunity cost.

We assume that the drivers are rational queue participants \citep{palm1957research}, who will join the queue when $U_{m}^n\geq 0$, and do not renege from the queue after joining it. Thus, when the queue length is $n$, the queue-joining process of group-$m$ drivers is also an independent Poisson process with rate $\tilde{\lambda}_m^n$,

\begin{equation}
\tilde{\lambda}_m^n=\lambda_m\times \mathbb{1}(U_{m}^n\geq 0)=\lambda_m\times \mathbb{1}\Bigl(R_m-p_m-r(T_d+W_m^n)\geq 0\Bigr), \quad \forall n\geq 0
\label{eq:join_rate}
\end{equation}

\subsection{The lottery control}

The platform controls the idle driver virtual queue by two levers. First, they can determine a set of commission fees $\{p_m\}$ for each group of drivers, which is independent of the status of the virtual queue to avoid dynamic pricing. This fee is collected from the drivers from every terminal-area order. Second, they can control a set of entry position lotteries $\{\Delta_m\}$ for each group of drivers that determine their starting position if they join the queue.

An entry position lottery $\Delta_m$ is a series of probability distributions $\{\delta_m^1,\delta_m^2,\cdots\}$ that works as a lottery to determine a group-$m$ driver's starting position in the queue. Each element $\delta_m^n$ in the set is a discrete probability distribution that specifies the lottery to be used when the queue length is $n$ upon a group-$m$ driver's arrival. 

\begin{figure}[!ht]
  \centering
  \includegraphics[width=0.9\textwidth]{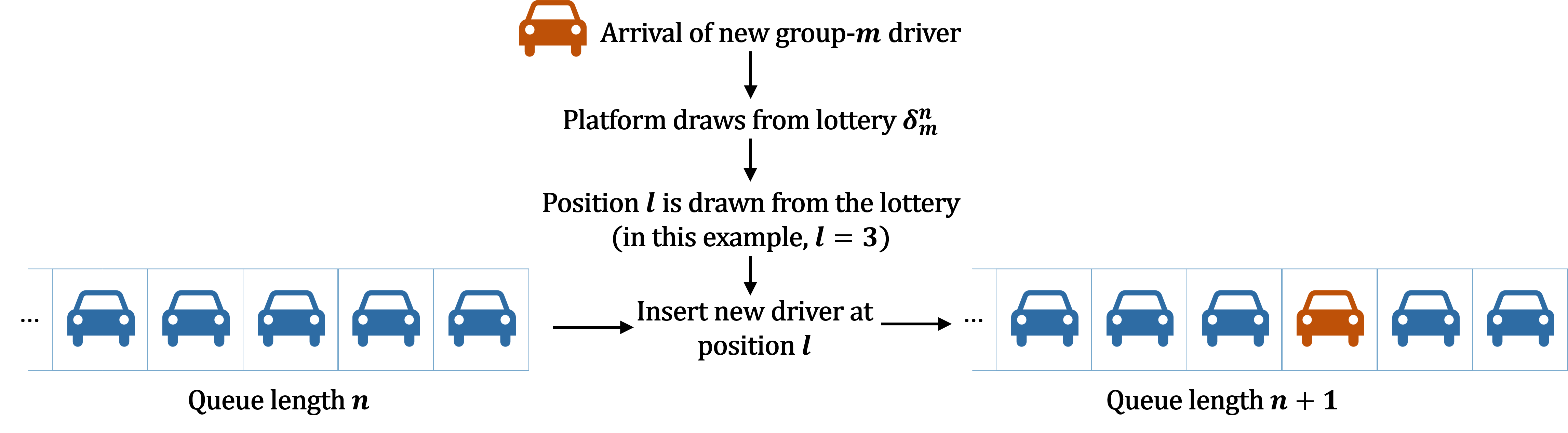}
  \caption{The entry position lottery system}\label{fig:lottery_system}
\end{figure}

When a group-$m$ driver arrives at the terminal, should the driver decide to join the virtual queue, the platform would draw a random integer $l\in \{1,2,\cdots,n+1\}$ using the lottery $\delta_m^n=\{\delta_m^{1,n},\delta_m^{2,n},\cdots,\delta_m^{n+1,n}\}$. In the lottery, the probability mass $\delta_m^{l,n}$ is the probability of the integer $l$ being drawn. After the draw, the driver will be inserted at the position $l$ in the queue. $l=1$ indicates that the driver is placed at the front of the queue, and $l=n+1$ indicates that the driver is placed at the back. At the same time, the drivers originally at queue positions $l,l+1,\cdots,n$ will each be moved one position back in the queue. An illustration of the operations of the lotteries is shown in \Cref{fig:lottery_system}.

By normalization, the elements of lottery $\delta_m^n$ satisfy,
\begin{equation}
    \sum_{l=1}^{n+1} \delta_m^{l,n}=1, \quad \forall n\geq 0
    \label{eq:distribution_normalization}
\end{equation}

\subsection{System dynamics and the steady state}
\begin{figure}[!ht]
  \centering
  \includegraphics[width=0.75\textwidth]{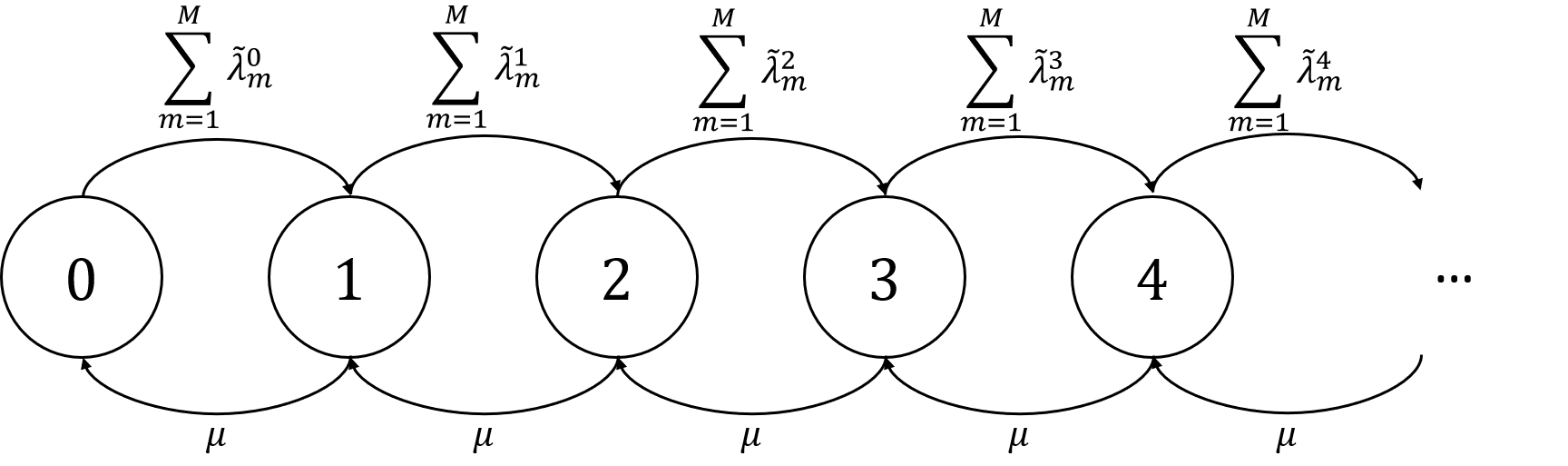}
  \caption{State transition diagram for the system dynamics}\label{fig:system_original_queue}
\end{figure}

Given the platform's $\{p_m\}$ and $\{\Delta_m\}$, the dynamics of the system can be illustrated by a queueing system with state-dependent arrival rates in \Cref{fig:system_original_queue}.

As given in \Cref{eq:join_rate}, the arrival rates of the queueing system in \Cref{fig:system_original_queue} are dictated by the drivers' joining behavior, which is controlled by the commission fees $\{p_m\}$ and the expected waiting times in the virtual queue $\{W_m^n\}$. The expected waiting times are impacted by the lotteries and the conditional expected waiting times $\{w_n^l\}$, which denotes the expected waiting time if a driver joins the queue at position $l$ when the queue length upon arrival is $n$. The relationship between $W_m^n$ and $w_n^l$ can be expressed by,
\begin{equation}
W_m^n=\sum_{l=1}^{n+1}\delta_m^{l,n}w_n^l \quad \forall m=1,2,\dots,M;\: n=1,2,\dots
\label{eq:expected_waiting_time}
\end{equation}

\begin{figure}[!ht]
  \centering
  \includegraphics[width=0.6\textwidth]{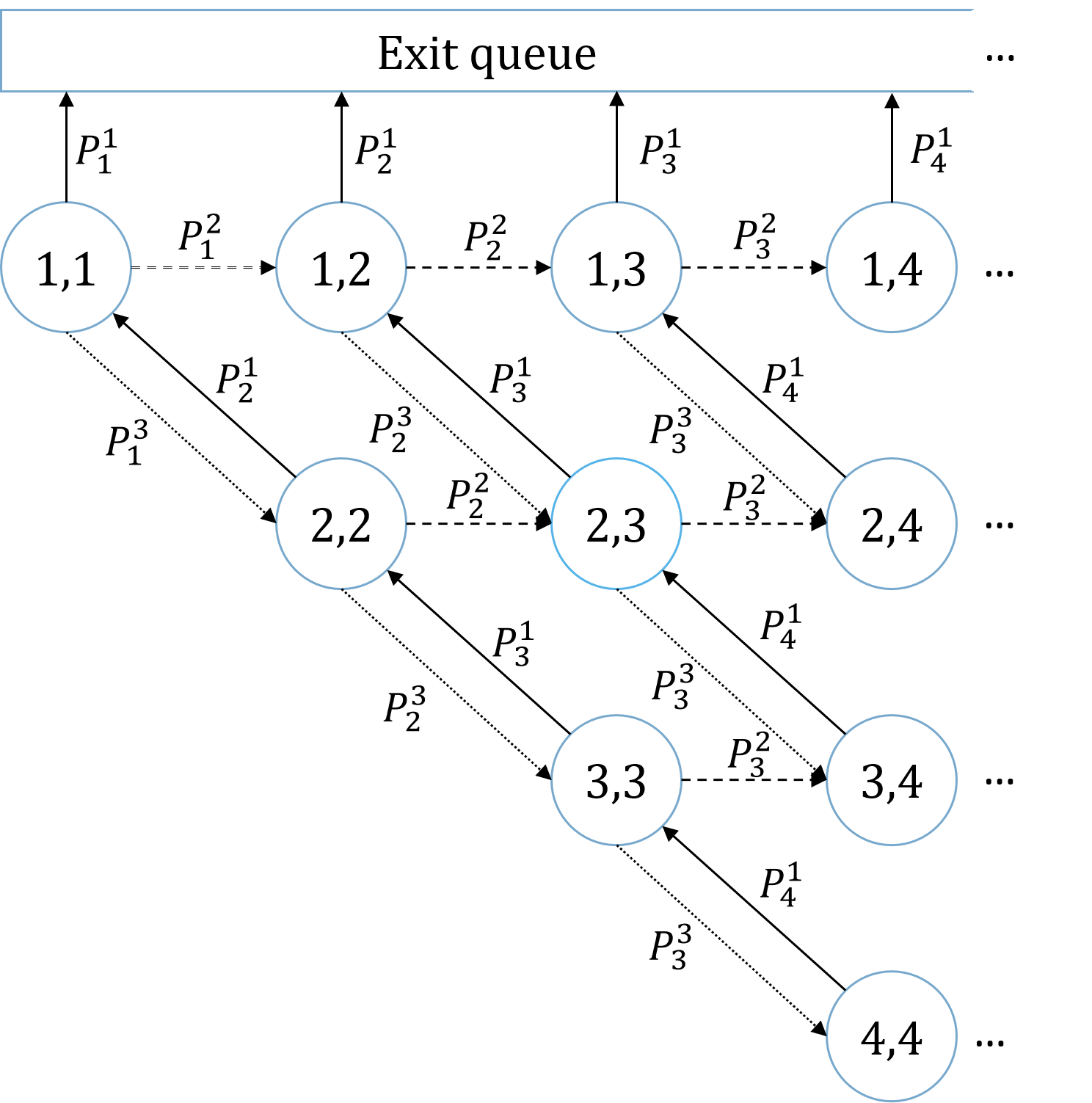}
  \caption{State transition probabilities for queue position and length}\label{fig:system_trans_queue_pos}
\end{figure}

To calculate the conditional waiting times $w_n^l$, we need to track the queueing progress of drivers by considering the subsequent arrivals and departures after a driver joins the queue. This progress can be tracked by the progression of a driver's queue position and the number of drivers in the queue, which can be represented by the transition diagram of states $(l,n)$ shown in \Cref{fig:system_trans_queue_pos}.

In \Cref{fig:system_trans_queue_pos}, the transition probabilities can be classified into three types:
\begin{itemize}
    \item $P_n^1$ is the probability of the next event in the queueing system being the arrival of a passenger. Therefore, if $l=1$, the driver gets matched and exits the queue. If $l>1$, the driver moves forward one position in the queue and the queue length reduces by one, making the next state $(l-1,n-1)$. The probability $P_n^1$ can be expressed by,
    \begin{equation*}
        P_n^1=\frac{\mu}{\sum_{c=1}^{M}\tilde{\lambda}_{c}^{n+1}+\mu}
    \end{equation*}
    \item $P_n^2$ is the probability of the next event in the queueing system being the arrival of a new driver and the new driver is placed behind the driver we are tracking. Therefore, the tracked driver stays in the same position $l$ while the queue length grows by one to $n+1$, making the next state $(l,n+1)$. The probability $P_n^2$ can be expressed as,
    \begin{equation*}
        P_n^2=\sum_{m=1}^{M}\left[\frac{\tilde{\lambda}_{m}^{n+1}}{\sum_{c=1}^{M}\tilde{\lambda}_{c}^{n+1}+\mu}(1-\sum_{i=1}^{l}\delta_m^{i,n+1})\right]
    \end{equation*}
    \item $P_n^3$ is the probability of the next event in the queueing system being the arrival of a new driver and the new driver is placed in front the driver we are tracking. Therefore, the tracked driver moves back one position to $l+1$ while the queue length grows by one to $n+1$, making the next state $(l+1,n+1)$. The probability $P_n^3$ can be expressed as,
    \begin{equation*}
        P_n^3=\sum_{m=1}^{M}\left(\frac{\tilde{\lambda}_{m}^{n+1}}{\sum_{c=1}^{M}\tilde{\lambda}_{c}^{n+1}+\mu}\sum_{i=1}^{l}\delta_m^{i,n+1}\right)
    \end{equation*}
\end{itemize}

Given the memory-less properties of the Poisson process, the conditional expected waiting time $w_n^l$ is the expected time to exit the queue for a driver in state $(l,n+1)$ in \Cref{fig:system_trans_queue_pos}. Therefore, by the recursive balance equation, the expected waiting time $w_n^l$ can be expressed as,

\begin{equation}
    \begin{aligned}
         w_n^l=&\sum_{m=1}^{M}\left[\frac{\tilde{\lambda}_{m}^{n+1}}{\sum_{c=1}^{M}\tilde{\lambda}_{c}^{n+1}+\mu}\left(\frac{1}{\sum_{c=1}^{M}\tilde{\lambda}_{c}^{n+1}+\mu}+\sum_{i=1}^{l}\delta_m^{i,n+1}w_{n+1}^{l+1}+(1-\sum_{i=1}^{l}\delta_m^{i,n+1})w_{n+1}^{l}\right)\right] &\\
        &\frac{\mu}{\sum_{c=1}^{M}\tilde{\lambda}_{c}^{n+1}+\mu}\times \left(\frac{1}{\sum_{c=1}^{M}\tilde{\lambda}_{c}^{n+1}+\mu}+w_{n-1}^{l-1}\right), \quad \forall n={1,2,3,\cdots},\; l\in\{2,\cdots,n+1\}
    \end{aligned}
\label{eq:expected_waiting_time_original_inner}
\end{equation}
\noindent in which $1/(\sum_{c=1}^{M}\tilde{\lambda}_{c}^{n+1}+\mu)$ is the expected time to the next event (arrival of a new driver or arrival of a new passenger); $w_{n+1}^{l+1}$ is the conditional expected waiting time after the new arrival is placed before the driver at position $l$; $w_{n+1}^{l}$ is the conditional expected waiting time after the new arrival is placed after the driver at position $l$; $w_{n-1}^{l-1}$ is the expected waiting time after the driver in front of the queue is matched to a passenger.

For $l=1$, the arrival of a passenger will cause the driver at the front to be matched. Therefore, \Cref{eq:expected_waiting_time_original_inner} reduces to,
\begin{equation}
\begin{aligned}
    w_n^1=&\sum_{m=1}^{M}\left[\frac{\tilde{\lambda}_{m}^{n+1}}{\sum_{c=1}^{M}\tilde{\lambda}_{c}^{n+1}+\mu}\left(\frac{1}{\sum_{c=1}^{M}\tilde{\lambda}_{c}^{n+1}+\mu}+\delta_m^{1,n+1}w_{n+1}^{2}+(1-\delta_m^{1,n+1})w_{n+1}^1\right)\right] \\
        &\frac{\mu}{\sum_{c=1}^{M}\tilde{\lambda}_{c}^{n+1}+\mu}\times \frac{1}{\sum_{c=1}^{M}\tilde{\lambda}_{c}^{n+1}+\mu}, \quad \forall n={0,1,2\cdots}
\end{aligned}
\label{eq:expected_waiting_time_original_boundary}
\end{equation}

Overall, \Cref{eq:join_rate,eq:distribution_normalization,eq:expected_waiting_time,eq:expected_waiting_time_original_inner,eq:expected_waiting_time_original_boundary} defines the equilibrium behavior of the virtual queue waiting times depicted in \Cref{fig:system_original_queue}. Under equilibrium, the steady-state probability $P_n$ of the queue length $n\in \{0,1,2,\cdots\}$ can be expressed as,
\begin{equation}
    \begin{aligned}
        &P_0=\frac{1}{1+\sum_{i=1}^{\infty}\left(\Pi_{j=1}^i \frac{\sum_{c=1}^M\tilde{\lambda}_{c}^{j-1}}{\mu}\right)}\\
        &P_n=\frac{\Pi_{j=1}^n \frac{\sum_{c=1}^M\tilde{\lambda}_{c}^{j-1}}{\mu}}{1+\sum_{i=1}^{\infty}\left(\Pi_{j=1}^i \frac{\sum_{c=1}^M\tilde{\lambda}_{c}^{j-1}}{\mu}\right)}, \quad \forall n=1,2,\cdots
    \end{aligned}
    \label{eq:steady_state_probability}
\end{equation}

\subsection{Formulation of the optimal control problem}

We consider the profit maximization problem and the social welfare optimization problem for the terminal ride-sourcing system. First, we derive the expected profit gain rate from the terminal area. The conditional profit gain rate of the system when the queue length is $n$ can be expressed as,
\begin{equation}
    \sum_{m=1}^M\tilde{\lambda}_{m}^n \bigl(p_m-\nu (T_d+W_m^n)\bigr)
    \label{eq:profit_rate_raw}
\end{equation}
\noindent in which $\nu$ is the profit rate for the platform in the general market outside of the terminal area, $p_m-\nu (T_d+W_m^n)$ is the platform's expected profit gain for a group-$m$ driver who joined the queue at length $n$, containing the commission fee and cost of driver idle and delivery time.

\begin{equation}
    \begin{aligned}
    \underset{\{p_m\},\{\Delta_m\}}{\max} \quad & \sum_{n=0}^{\infty} P_n \sum_{m=1}^M \Bigl(\tilde{\lambda}_{m}^n(p_m-\nu (T_d+W_m^n))\Bigr)  \\
    \text{subject to} \quad & (\ref{eq:join_rate}), (\ref{eq:distribution_normalization}), (\ref{eq:expected_waiting_time}), (\ref{eq:expected_waiting_time_original_inner}), (\ref{eq:expected_waiting_time_original_boundary}),(\ref{eq:steady_state_probability})
    \label{eq:profit_optimization}
    \end{aligned}
\end{equation}
\noindent which is an infinite-dimensional non-linear optimization problem. 

For social welfare maximization, similar to \Cref{eq:profit_rate_raw}, the social-welfare gain rate conditional on queue length $n$ is:
\begin{equation}
    \sum_{m=1}^M\tilde{\lambda}_{m}^n(R_m-(r+\nu) (T_d+W_m^n))
\end{equation}
\noindent in which $R_m-(\nu+r) (T_d+W_m^n)$ is the additional expected social welfare gained by a group-$m$ driver joining the queue at length $n$.

Consequently, the social welfare optimization problem can be formulated as,
\begin{equation}
\begin{aligned}
    & \underset{\{p_m\},\{\Delta_m\}}{\max}
    && \sum_{n=0}^{\infty} P_n\sum_{m=1}^M\Bigl(\tilde{\lambda}_{m}^n(R_m-(r+\nu) (T_d+W_m^n))\Bigr) \\
     &\text{subject to}  &&(\ref{eq:join_rate}), (\ref{eq:distribution_normalization}), (\ref{eq:expected_waiting_time}), (\ref{eq:expected_waiting_time_original_inner}), (\ref{eq:expected_waiting_time_original_boundary}), (\ref{eq:steady_state_probability})
\end{aligned}
\label{eq:welfare_optimization}
\end{equation}
\noindent which, similar to problem (\ref{eq:profit_optimization}), is also an infinite-dimensional non-linear optimization problem.


\section{Equilibrium Property Analysis and Problem Reformulation}
\label{sec:theory}

While \Cref{sec:model} defines the system model and optimal control problems (\ref{eq:profit_optimization}) and (\ref{eq:welfare_optimization}), both the system equilibrium and the optimization problems are challenging to solve. The main challenges here are twofold:

\begin{enumerate}
    \item The queueing system in \Cref{fig:system_original_queue} has state-dependent arrival rates and infinite capacity. Therefore, the lottery $\{\Delta_m\}$ is a variable of infinite dimensions.
    \item The equilibrium defined by \Cref{eq:join_rate,eq:distribution_normalization,eq:expected_waiting_time,eq:expected_waiting_time_original_inner,eq:expected_waiting_time_original_boundary} involves non-smooth and nonlinear functions, which increases the difficulty of solving for the equilibrium.
\end{enumerate}

In this section, we aim to address the two main challenges by analyzing the structural properties of the optimal equilibrium. We start by decomposing the structure of the optimal equilibrium into two parts in \Cref{subsec:theo_cha_1} and transform the optimal control problems into bi-level programming problems. Then, we study the structural properties of equilibrium arrival rates in \Cref{subsec:theo_cha_2} to further simplify the candidate equilibrium set. Taking these insights, we present the transformed optimal control problems and the solution methodology in \Cref{subsec:theo_cha_3}, followed by a theoretical analysis of the performance of the lottery control and our problem transformation in \Cref{subsec:theo_cha_4}.

\subsection{Decomposing the optimal equilibrium} \label{subsec:theo_cha_1}
First, we address the challenge of the infinite dimensionality of the equilibrium. Due to the recursive nature of \Cref{eq:expected_waiting_time_original_inner} and \Cref{eq:expected_waiting_time_original_boundary}, to solve the equilibrium of the queueing system, the lottery $\{\Delta_m\}$ must be fully specified.  Therefore, we take an alternative approach and first analyze the structure of the optimal equilibrium, and then explore the lotteries that can generate such equilibria to decompose the problem and reduce its size.

First, we consider the equilibrium that can maximize the platform's profit defined in the optimization problem (\ref{eq:profit_optimization}). We first show that a sufficient condition on the equilibrium queue-joining rate $\{\tilde{\lambda}_m^n\}$ exists for the optimal solution in \Cref{prop:threshold_profit_max}. 

\begin{prop}
     The equilibrium queue-joining rate $\{\tilde{\lambda}_m^n\}$ of the profit-maximizing solution of problem (\ref{eq:profit_optimization}) satisfies,
    \begin{center}
$\displaystyle \exists N_{m,p}^* \geq 1 \text{ such that } \forall n \geq N_{m,p}^*, \tilde{\lambda}_m^n = 0.$
\end{center}
\label{prop:threshold_profit_max}
\end{prop}

\begin{figure}[!ht]
  \centering
  \includegraphics[width=0.75\textwidth]{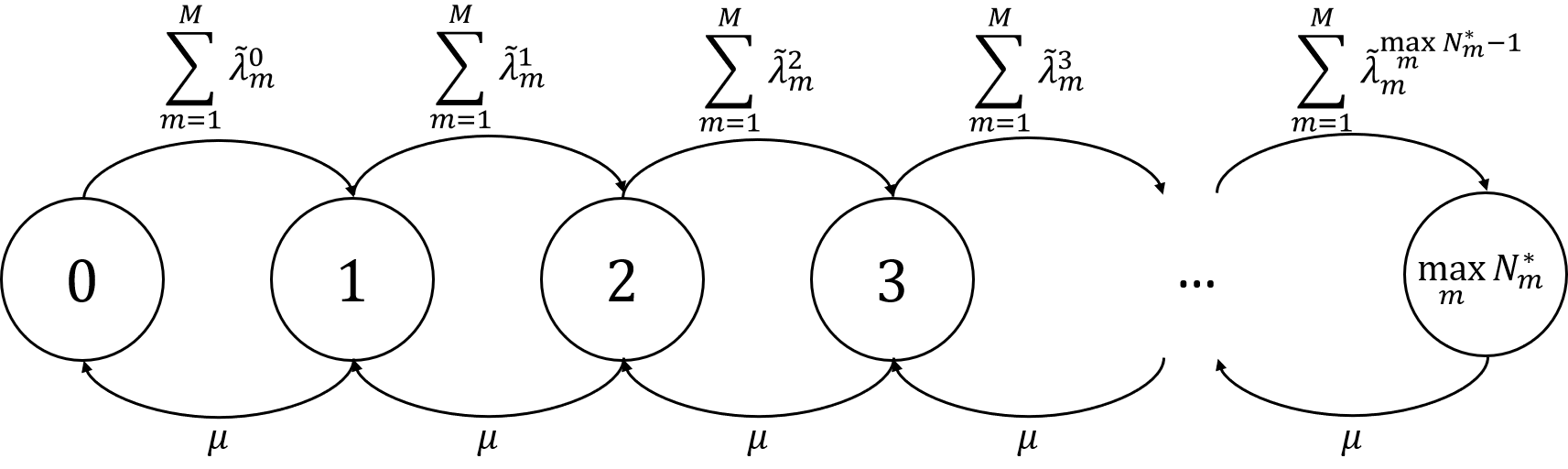}
  \caption{State transition diagram for the profit-maximizing equilibrium}\label{fig:system_reduced_queue}
\end{figure}

\Cref{prop:threshold_profit_max} states that the profit-maximizing system equilibrium will have a limited capacity. Therefore, the queue structure of the optimal solution can be reduced to the queue shown in \Cref{fig:system_reduced_queue}.

To optimize the system profit, the platform's strategy aligns with the equilibrium point of the system in \Cref{fig:system_reduced_queue}. As the queue lengthens, the platform's opportunity cost eventually outweighs the commission fee, making further admissions unprofitable. The net profit from an additional driver is shown to be a decreasing function of the queue length. As the queue grows, this marginal profit diminishes and eventually drops below zero. This framework helps identify the optimal operational capacity for each driver group, with an upper bound presented in \Cref{cor:bound_profit_max}.

\begin{cor}
For any given price $p_m$, the optimal profit-maximizing threshold $N_{m,p}^*$ has the following upper bound,
    \begin{center}
    $\displaystyle N_{m,p}^*\leq \Bigl\lceil \frac{\mu (p_m-vT_d)}{\nu}-\frac{1}{2}\Bigr\rceil$
\end{center}
    \label{cor:bound_profit_max}
\end{cor}

In \Cref{cor:bound_profit_max}, we note that the upper bound of the capacity $N_{m,p}^*$ is positively correlated to the customer arrival rate $\mu$ and commission fee $p_m$ and negatively correlated with the profit rate for the platform in the general market $\nu$. This indicates that as customers arrive more often, the waiting time of idle drivers is reduced, and the platform can afford to allow more drivers into the queue. Similarly, an increment in commission fee increases the platform's profit margin and leads to the platform allowing more idle drivers in the queue. An increment in the profitability of the general market means that the idle profit loss rate of the platform increases, which leads to the platform retaining less idle drivers at the terminal, as the cost of retaining each driver would be higher.

Similarly, on the social-welfare-maximizing equilibrium solution defined in problem (\ref{eq:welfare_optimization}), we can also obtain a sufficient condition of the equilibrium queue-joining rate in \Cref{prop:threshold_welfare_max}.

\begin{prop}
    The equilibrium queue-joining rate $\{\tilde{\lambda}_m^n\}$ of the social-welfare-maximizing solution of problem (\ref{eq:welfare_optimization}) satisfies,
    \begin{center}
$\displaystyle \exists N_{m,s}^* \geq 1 \text{ such that }\forall n \geq N_{m,s}^*, \tilde{\lambda}_m^n = 0$
\end{center}
\label{prop:threshold_welfare_max}
\end{prop}

\Cref{prop:threshold_welfare_max} can be interpreted in similar ways as \Cref{prop:threshold_profit_max}. We can also obtain an upper bound on the capacity $N_{m,s}^*$ for the social welfare optimization problem in \Cref{cor:bound_welfare_max}.

\begin{cor}
   For any given price $p_m$, the optimal profit-maximizing threshold $N_{m,s}^*$ has the following upper bound,
   \begin{center}
    $\displaystyle N_{m,s}^*\leq \Bigl\lceil \frac{\mu (R_m-(r+v)T_d)}{\nu+r}-\frac{1}{2}\Bigr\rceil$
\end{center}
    \label{cor:bound_welfare_max}
\end{cor}

Next, we discuss the design of the decision variables $\{p_m\}$ and $\{\Delta_m\}$ that leads to the structure of the equilibrium joining behavior presented in \Cref{prop:threshold_profit_max} and \Cref{prop:threshold_welfare_max}. To achieve this, the decision variables would have to ensure $U_{m}^n<0$ for all $n\geq N_{m,s}^*$. However, since the commission fee does not change with respect to queue length, we can only control the lotteries $\{\Delta_m\}$ to achieve this purpose. By \Cref{eq:join_rate}, we need to design the lotteries such that the expected waiting time $W_m^n$ is sufficiently long when the queue length exceeds the optimal capacity to induce balking. In \Cref{prop:extension_policy}, we prove that for any capacity 
$N$, it is feasible to design a lottery such that the waiting time is longer for anyone joining the queue at or beyond its capacity, compared to someone joining when the queue is below capacity, thereby ensuring the system attains the desired equilibrium state.

\begin{prop}
    For any given $\{\delta_m^1,\delta_m^2,\cdots,\delta_m^N\}$, 
    $\displaystyle \exists \{\delta_m^{N+1}, \delta_m^{N+2}, \cdots\}$ \text{ such that } \begin{center}
    $\forall n \geq N$,\; $ W_m^n > \max_{i<N}(W_m^i)$    
    \end{center}
            
\label{prop:extension_policy}
\end{prop}

\Cref{prop:extension_policy} introduces two key insights. It demonstrates that when the queue is at capacity, wait times can be managed to avoid system overload. Therefore, by \Cref{eq:join_rate}, once the optimal capacity $N_m^*$ is fixed, we can coordinate the price $p_m$ with the lottery $\Delta_m$ such that,
\begin{equation}
    R_m-p_m-r(T_d+\max_{n<N_m^*}W_m^n)= 0
    \label{eq:capa_zero_utility}
\end{equation}
\noindent which ensures all group-$m$ drivers will balk from the queue when $n\geq N_m^*$. Furthermore, \Cref{prop:extension_policy} allows us to decouple the lottery design of those under-capacity lottery items $\{\delta_m^1,\delta_m^2,\cdots,\delta_m^N\}$ and over-capacity items $\{\delta_m^{N+1},\delta_m^{N+2},\cdots,\}$. Combining \Cref{prop:threshold_profit_max} or \Cref{prop:threshold_welfare_max} and \Cref{prop:extension_policy}, the original optimization problems \ref{eq:profit_optimization} and \ref{eq:welfare_optimization} can be transformed to a bi-level problem in which,
\begin{itemize}
    \item On the upper level, we search for the optimal capacities.
    \item In the lower-level problem, given the capacities determined by the upper level, we can directly solve for the optimal price and the corresponding under-capacity lottery policies, satisfying the condition in \Cref{eq:capa_zero_utility}. The over-capacity policies, $\{\delta_m^{N_m}, \delta_m^{N_m+1}, \cdots\}$, can be found according to \Cref{prop:extension_policy} to complete the overall solution.
\end{itemize}

This allows us to transform the original problems into a bi-level optimization problem, in which the equilibrium solution is decomposed into two segments and only the finite-dimensional segment needs to be solved. Therefore, our first main challenge is addressed.

\subsection{Structural properties of arrival rates} \label{subsec:theo_cha_2}

However, solving for the queueing system's equilibrium is still challenging due to the non-linear and non-smooth characteristics of the governing equations. We overcome this difficulty by analyzing the structure of the optimal equilibrium to constrain the feasible set and eliminate these non-smooth characteristics effectively.

A key insight from the proofs of \Cref{prop:threshold_profit_max} and \Cref{prop:threshold_welfare_max} is that the marginal profit or welfare from admitting a new driver decreases as the queue length increases. This implies that with direct control over admissions, the optimal policy would be to accept all drivers in a group as long as the queue length is below its capacity, $N_m^*$. In our model, however, the platform's control is limited to setting prices and lotteries, which cannot guarantee this ideal outcome. Nevertheless, as shown in \Cref{prop:threshold_all_join_profit}, some profit-maximizing solutions do adhere to this structure under specific price conditions.

\begin{prop}
    For any given price $p_m\leq R_m-rT_d-\frac{r}{\mu-\sum_{c=1}^M\lambda_c}$, the profit-maximizing equilibrium arrival rates $\tilde{\lambda}_m^n$ satisfy,\\
    \begin{center}
    $\displaystyle \exists N_{m,p}^* \text{ such that } \tilde{\lambda}_m^n =
    \begin{cases}
        \lambda_m & \text{if } n < N_{m,p}^* \\
        0       & \text{otherwise}
    \end{cases}$
\end{center}
\label{prop:threshold_all_join_profit}
\end{prop}

\Cref{prop:threshold_all_join_profit} shows that for some prices, the optimal lottery conditional on the prices will result in all group-$m$ drivers joining the queue when the queue length is under the optimal capacity $N_{m,p}^*$. Similarly, for the social welfare optimization problem, we can obtain a similar result in \Cref{prop:threshold_all_join_welfare},

\begin{prop}
    For any given price $p_m\leq R_m-rT_d-\frac{r}{\mu-\sum_{c=1}^M\lambda_c}$, the social-welfare-maximizing equilibrium arrival rates $\tilde{\lambda}_m^n$ satisfy,\\
    \begin{center}
    $\displaystyle \exists N_{m,s}^* \text{ such that } \tilde{\lambda}_m^n =
    \begin{cases}
        \lambda_m & \text{if } n < N_{m,s}^* \\
        0       & \text{otherwise}
    \end{cases}$
\end{center}
\label{prop:threshold_all_join_welfare}
\end{prop}

Both \Cref{prop:threshold_all_join_profit} and \Cref{prop:threshold_all_join_welfare} suggest that, at least some price-conditional optimal solutions have a threshold structure, where all drivers join below a certain queue length and none join above it. Taking this insight, instead of solving the original optimization problems, if we further reduce the feasible region to constrain the resulting equilibrium to have such a structure, the solution may still have high quality. Furthermore, this approach can address the solution challenges regarding the lower-level equilibrium in the aforementioned bi-level programming approach. When the capacities $N_m$ are fixed by the upper level, imposing this threshold structure simplifies the queue-joining rate to:
\begin{equation}
    \tilde{\lambda}_m^n=\lambda_m \mathbb{1}(n<N_m)
    \label{eq:conditional_arrival_rate_new}
\end{equation}
\noindent which make the queue joining rate $\tilde{\lambda}_m^n$ becoming deterministic after $N_m$s are fixed in the upper level. This formulation eliminates the non-smoothness present in \Cref{eq:join_rate}, making the equilibrium significantly more tractable.

\subsection{Problem transformation and solution methodology} \label{subsec:theo_cha_3}
Given the insights in \Cref{subsec:theo_cha_1} and \Cref{subsec:theo_cha_2}. The structural constraint on the arrival rates in \Cref{subsec:theo_cha_2} changes the conditional waiting time functions to,
\begin{subequations}
\label{eq:conditional_waiting_time_bounded}
\begin{align}
    w_n^l=&\sum_{m=1}^{M}\left[\frac{\tilde{\lambda}_{m}^{n+1}}{\sum_{c=1}^{M}\tilde{\lambda}_{c}^{n+1}+\mu}\left(\frac{1}{\sum_{c=1}^{M}\tilde{\lambda}_{c}^{n+1}+\mu}+\sum_{i=1}^{l}\delta_m^{i,n+1}w_{n+1}^{l+1}+(1-\sum_{i=1}^{l}\delta_m^{i,n+1})w_{n+1}^{l}\right)\right]  \label{eq:conditional_waiting_time_bounded_regular_under_capacity}\\
    &+\frac{\mu}{\sum_{c=1}^{M}\tilde{\lambda}_{c}^{n+1}+\mu}\times \left(\frac{1}{\sum_{c=1}^{M}\tilde{\lambda}_{c}^{n+1}+\mu}+w_{n-1}^{l-1}\right), \quad  \forall n={1,2,\cdots,\hat{N}-2},\; l\in\{2,3,\cdots,n+1\} \notag \\
     w_n^1=&\sum_{m=1}^{M}\left[\frac{\tilde{\lambda}_{m}^{n+1}}{\sum_{c=1}^{M}\tilde{\lambda}_{c}^{n+1}+\mu}\left(\frac{1}{\sum_{c=1}^{M}\tilde{\lambda}_{c}^{n+1}+\mu}+\delta_m^{1,n+1}w_{n+1}^{2}+(1-\delta_m^{1,n+1})w_{n+1}^1\right)\right]  \label{eq:conditional_waiting_time_bounded_first_under_capacity}\\
        &\frac{\mu}{\sum_{c=1}^{M}\tilde{\lambda}_{c}^{n+1}+\mu}\times \frac{1}{\sum_{c=1}^{M}\tilde{\lambda}_{c}^{n+1}+\mu}, \quad \forall n={0,1,2,\cdots,\hat{N}-2} \notag\\
    w_{\hat{N}-1}^l&=\frac{1}{\mu}+w_{\hat{N}-2}^{l-1}, \quad  \forall l\in {2,3,\cdots,\hat{N}} \label{eq:conditional_waiting_time_bounded_regular_at_capacity}\\
        w_{\hat{N}-1}^1&=\frac{1}{\mu} \label{eq:conditional_waiting_time_bounded_first_at_capacity}
\end{align}
\end{subequations}
\noindent of which $\hat{N}=\max_mN_m$ is the rejection threshold of the queue, \Cref{eq:conditional_waiting_time_bounded_regular_under_capacity} and \Cref{eq:conditional_waiting_time_bounded_first_under_capacity} follows \Cref{eq:expected_waiting_time_original_inner} and \Cref{eq:expected_waiting_time_original_boundary}, correspondingly, to capture the expected conditional waiting time of drivers at different queue positions when the queue is under capacity. When the queue is at capacity after the driver joins, no new arrivals will join the queue. Therefore, if the driver joins at first position, the driver is guaranteed to be matched when the next rider arrives, which is captured by \Cref{eq:conditional_waiting_time_bounded_first_at_capacity}. Otherwise, the driver will advance one position at the next arrival, which is modeled by \Cref{eq:conditional_waiting_time_bounded_regular_at_capacity}.

For profit maximization, we can transform the original optimization problem (\ref{eq:profit_optimization}) to an alternative bi-level program (\ref{eq:profit_maximization_new}),
\begin{equation}
\begin{aligned}
    & \max_{\{N_m\}} 
    && \sum_{m=1}^M\lambda_m\sum_{n=1}^{N_m}\left(R_m-(r+\nu )T_d-r\xi_m^*-\nu W_m^{n*}\right)\\
    & \text{s.t.} 
    && N_m \in \mathbb{Z}^+, \quad \forall m\in\{1,2,\cdots,M\}\\ 
    & \text{where}\; \Delta_m^*,\xi_m^*,W_m^{n*} = 
    && \argmax_{\{p_m\},\{\Delta_m\},\{\xi_m\},\{W_m^n\}}\sum_{m=1}^M\lambda_m\sum_{n=1}^{N_m}\left(R_m-(r+\nu )T_d-r\xi_m-\nu W_m^n\right)\\
    & 
    && \quad \text{s.t.} \quad (\ref{eq:distribution_normalization}),(\ref{eq:expected_waiting_time}),(\ref{eq:steady_state_probability}), (\ref{eq:conditional_arrival_rate_new}),(\ref{eq:conditional_waiting_time_bounded})\\
    & 
    && \quad W_m^n\leq \xi_m, \quad \forall m\in\{1,2,\cdots,M\},\; n\in\{0,1,\cdots,N_m-1\}
\end{aligned}
\label{eq:profit_maximization_new}
\end{equation}

\noindent In program \eqref{eq:profit_maximization_new}, we reformulate the problem by replacing the commission fee variable $p_m$ with a new decision variable, $\xi_m$, which represents the maximum tolerable expected waiting time for drivers in group $m$. This is because in the profit maximization case, the platform sets the highest possible fee, which is determined by the maximum wait time drivers will accept, given by the equation $p_m = R_m - r(T_d + \xi_m)$. We use this relationship to substitute $p_m$ and add a constraint $W_m^n\leq \xi_m$ to ensure the wait time remains within this maximum limit.

Similarly, the social welfare maximization problem \ref{eq:welfare_optimization} can be reformulated as problem (\ref{eq:welfare_maximization_new}),
\begin{minipage}{\linewidth}
\begin{align}
    & \max_{\{N_m\}} 
    && \sum_{m=1}^M\lambda_m\sum_{n=1}^{N_m}\Bigl(R_m-(r+\nu )(T_d+W_m^{n*})\Bigr) \notag \\
    & \text{s.t.} 
    && N_m \in \mathbb{Z}^+, \quad \forall m\in\{1,2,...,M\} \notag \\
    & \text{where}\; p_m^*,\Delta_m^*,W_m^{n*} = 
    && \argmax_{\{p_m\},\{\Delta_m\},\{W_m^n\}}\sum_{m=1}^M\lambda_m\sum_{n=1}^{N_m}\Bigl(R_m-(r+\nu )(T_d+W_m^n)\Bigr) \label{eq:welfare_maximization_new} \\
    & 
    && \quad \text{s.t.} \quad (\ref{eq:distribution_normalization}), (\ref{eq:expected_waiting_time}), (\ref{eq:steady_state_probability}), (\ref{eq:conditional_arrival_rate_new}), (\ref{eq:conditional_waiting_time_bounded}) \notag \\
    & 
    && \quad R_m-p_m-r(T_d+W_m^n)\geq 0, \quad \forall m\in\{1,2,...,M\}, n\in\{0,1,...,N_m-1\} \notag\\
    \vspace{0.5em}\notag
\end{align}
\end{minipage}
\noindent in which the additional constraint $R_m-p_m-r(T_d+W_m^n)\geq 0$ ensures that all group-$m$ drivers would join the queue under the capacity $N_m$.

We use a genetic algorithm described in \Cref{alg:ga} to solve (\ref{eq:profit_maximization_new}) and (\ref{eq:welfare_maximization_new}).

\begin{algorithm}[h!]
\caption{Genetic Algorithm for solving (\ref{eq:profit_maximization_new}) and (\ref{eq:welfare_maximization_new})}
\label{alg:ga}
\begin{algorithmic}[1]
    \STATE \textbf{Initialize} population $\mathcal{P}_0 = \{\mathbf{N}_1, \dots, \mathbf{N}_P\}$ with random chromosomes.
    \FOR{$g = 1$ \TO $G_{\text{max}}$}
        \STATE \textbf{Evaluate fitness} $f(\mathbf{N})$ for each chromosome $\mathbf{N} \in \mathcal{P}_{g-1}$. By solve lower-level problem and use the optimum as the fitness measure.
        \STATE \textbf{Select} elite set $\mathcal{E} \gets \text{Top } E \text{ solutions from } \mathcal{P}_{g-1}$.
        \STATE \textbf{Initialize} next generation $\mathcal{P}_g \gets \mathcal{E}$.

        \WHILE{$|\mathcal{P}_g| < P$}
            \STATE Select parents $\mathbf{N}^A, \mathbf{N}^B \gets \text{TournamentSelect}(\mathcal{P}_{g-1})$.
            \STATE Generate a new child chromosome $\mathbf{N}^{\text{child}}$.
            \FOR{$m = 1$ \TO $M$}
                \STATE $N_m^{\text{child}} \gets \text{RandomChoice}(N_m^A, N_m^B)$.
                \IF{\text{rand()} $<$ $p_{\text{mut}}$}
                    \STATE Draw a random integer $\delta$ from $\{-1, 0, 1\}$.
                    \STATE $N_m^{\text{child}} \gets \max(1, N_m^{\text{child}} + \delta)$.
                \ENDIF
            \ENDFOR
        \ENDWHILE
        \STATE Add $\mathbf{N}^{\text{child}}$ to $\mathcal{P}_g$.
        \IF{$\max_{\mathbf{N}} f(\mathbf{N})$ is unchanged for $K$ generations} 
            \STATE \textbf{break}
        \ENDIF
    \ENDFOR
   \STATE \textbf{Return} $\mathbf{N}^* \gets \arg\max_{\mathbf{N} \in \mathcal{P}_g} f(\mathbf{N})$.
\end{algorithmic}
\end{algorithm}

In \Cref{alg:ga}, $\mathbf{N}_i$ is a $M$-dimension integer vector that contains the rejection capacity $N_1,N_2,\cdots,N_M$ for all driver groups; $P$ is the population size; \text{TournamentSelect()} is a selection process where a random subset of chromosomes is chosen from the population and the fittest individual from that subset is selected; \text{RandomChoice()} is a uniform crossover operator that constructs a child chromosome by randomly choosing each of its gene values from one of the two parent chromosomes; and $K$ is a hyperparameter that triggers the algorithm's convergence check.

To initialize the genetic permutation's efficiency, we first determine a set of upper bound of capacities $\{N_m^h\}$ solely for initialization. First, we perform a loose grid search by iteratively increasing each $N_m$ in coarse increments for the upper level. For each grid point, the lower-level problem is solved to evaluate the objective function. The search terminates when the marginal change of the profit turns negative, and we use that point as each $N_m^h$ for group-$m$ drivers.


\subsection{Performance of the problem transformation and lottery control} \label{subsec:theo_cha_4}

Problems (\ref{eq:profit_maximization_new}) and (\ref{eq:welfare_maximization_new}) are reduced form of problems (\ref{eq:profit_optimization}) and (\ref{eq:welfare_optimization}), correspondingly. However, in the case of social welfare maximization, we prove in \Cref{threorem:welfare_equivalence} that the restrictions in (\ref{eq:welfare_maximization_new}) do not prevent the attainment of the global optimum.
\begin{prop}
    The optimal solutions to problems (\ref{eq:welfare_maximization_new}) and (\ref{eq:welfare_optimization}) coincide.
    \label{threorem:welfare_equivalence}
\end{prop}
\Cref{threorem:welfare_equivalence} is significant because it demonstrates that the simplifications introduced in (\ref{eq:welfare_maximization_new}) do not compromise optimality, ensuring that the problem remains computationally efficient while still achieving the best possible social welfare. This equivalence validates the use of the reduced formulation without loss of generality, facilitating analytical and numerical tractability.

The performance of the lottery control is benchmarked against the FIFO dynamic pricing method \citep{chen2001state}. In this FIFO scheme, the platform sets a commission fee $p_m^n$ based on the queue length $n$ and places all arriving drivers at the back of the queue. \Cref{prop:welfare_lot_vs_price} demonstrates that, in comparison, lottery control yields an equivalent or superior social welfare outcome,
\begin{prop}
    The proposed approach (problem (\ref{eq:welfare_maximization_new})) yields an optimal social welfare greater than or equal to that of the dynamic pricing strategy.
    \label{prop:welfare_lot_vs_price}
\end{prop}

\Cref{prop:welfare_lot_vs_price} demonstrates the competitive advantage of the lottery control mechanism. The intuition behind this result is that lottery control provides more granular control over driver waiting times, which is the primary source of welfare loss in problem (\ref{eq:welfare_maximization_new}). Under the FIFO dynamic pricing policy, these waiting times are fixed once system capacities are set. In contrast, the lottery mechanism allows waiting times to be fine-tuned independently of capacity. Consequently, lottery control can achieve higher social welfare than the FIFO approach by more effectively minimizing driver idleness.


\section{Numerical Example} \label{sec:numerical}
In this section, we use numerical examples to demonstrate the solutions, properties, and performance of the lottery control. The operational setting we consider here is the operation of an idle driver queue at a Chinese metropolitan airport. At the airport area, the platform maintains a ``digital fence" and only drivers who are idle in the fence area can join the queue and be eligible to be matched to an arriving airport passenger. The arriving drivers are classified into two categories: involuntary drivers, who arrive at the airport area carrying a passenger; and voluntary drivers, who arrive at the airport area idle. The latter values the airport order higher due to the sunk cost associated with their idle repositioning process to the airport \citep{liu2025bounded}. We use $m=1$ and $m=2$ to denote the two types of drivers, respectively. We set up the model parameters using an operational dataset from this airport, which is described in detail by \cite{liu2025bounded}.

\begin{table}[h]\centering
\caption{Value of model parameters}
\label{tb:para_model}\footnotesize
\begin{tabular}{c c}
\hline
Parameter & Value\\
\hline
$\lambda_1$ & 31.3 veh/h\\
$\lambda_2$  & 10.6 veh/h\\
$\mu$ & 46.1 orders/h\\
$R_1$  & 80 CNY\\
$R_2$  & 90 CNY\\
$r$  & 40 CNY/h\\
$v$  & 5 CNY/h\\
$T_d$ & 30 minutes\\
\hline
\end{tabular}
\end{table}

\subsection{Solution Properties}
First, we showcase the profit-maximizing and the social-welfare-maximizing lotteries in the market described by \Cref{tb:para_model} to demonstrate the structures and properties of the optimal solutions and resulting equilibria.

\begin{figure}[h!]
\centering
    \subfloat[][Involuntary drivers]{\includegraphics[width=0.5\linewidth]{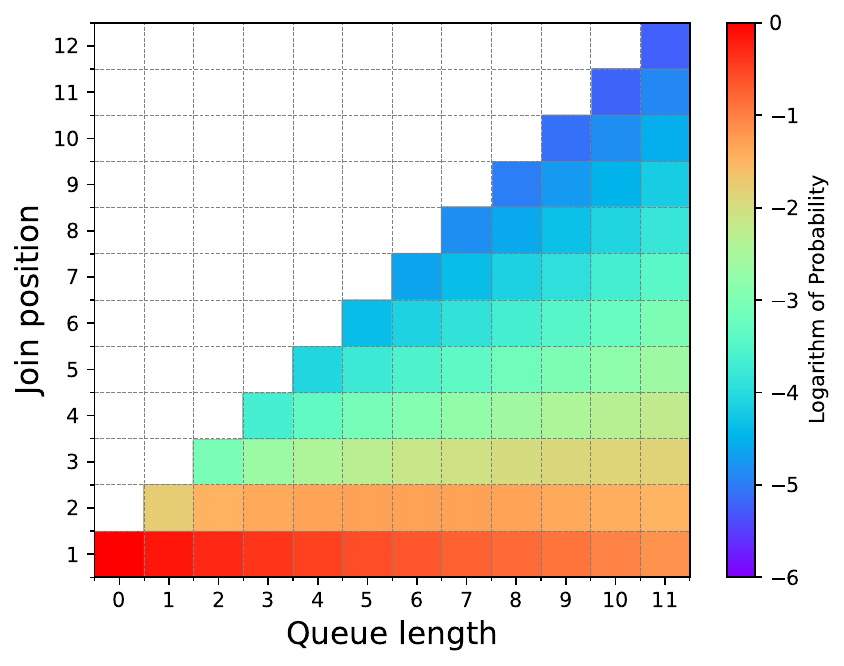}}
	\subfloat[][Voluntary drivers]{\includegraphics[width=0.5\linewidth]{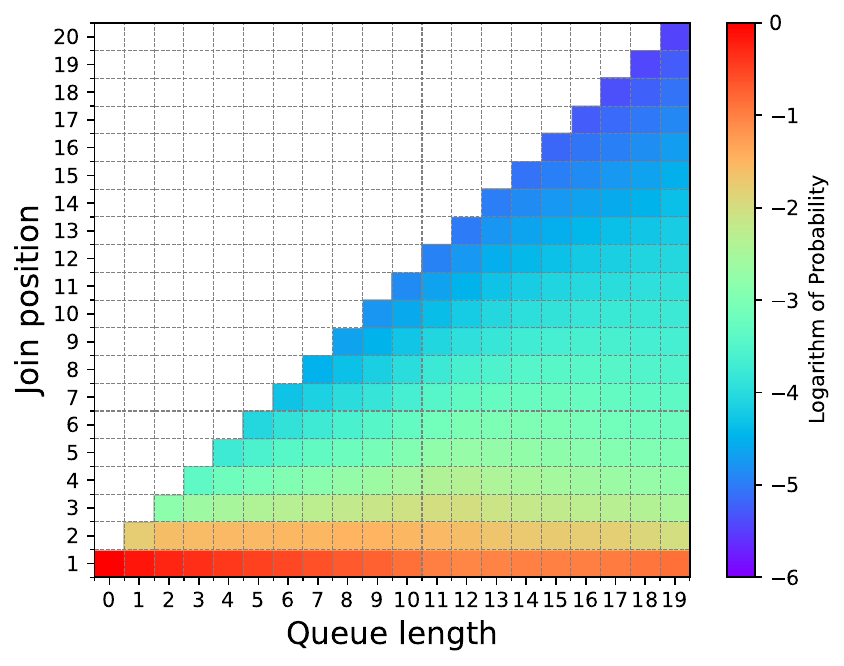}}
	\caption{The profit-maximizing driver entry position lotteries}
	\label{fig:numer_MO_lotteries}
\end{figure}

\begin{figure}[h!]
  \centering
  \includegraphics[width=0.6\textwidth]{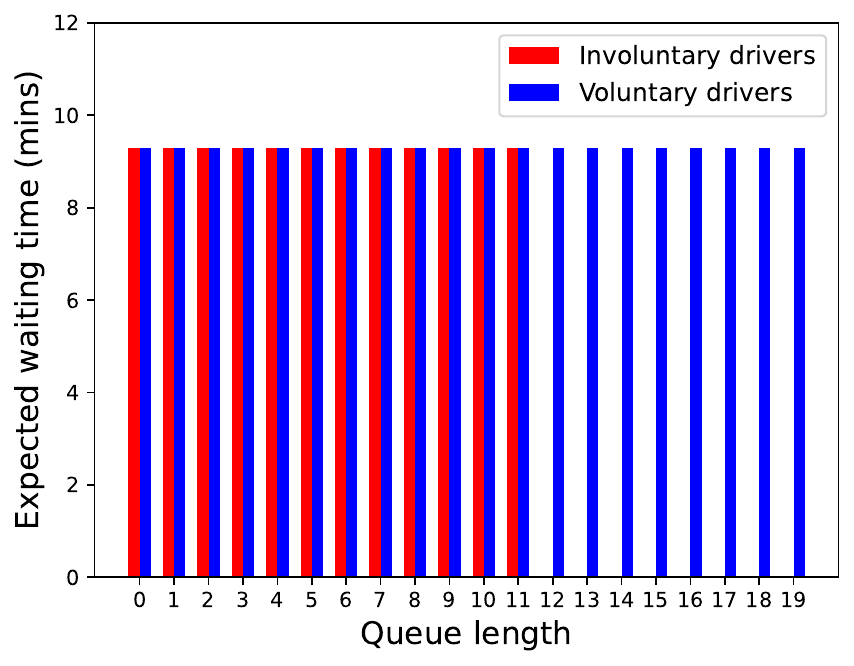}
  \caption{Equilibrium conditional waiting times for profit maximization}\label{fig:numer_PO_waittime}
\end{figure}

\begin{figure}[h!]
  \centering
  \includegraphics[width=0.6\textwidth]{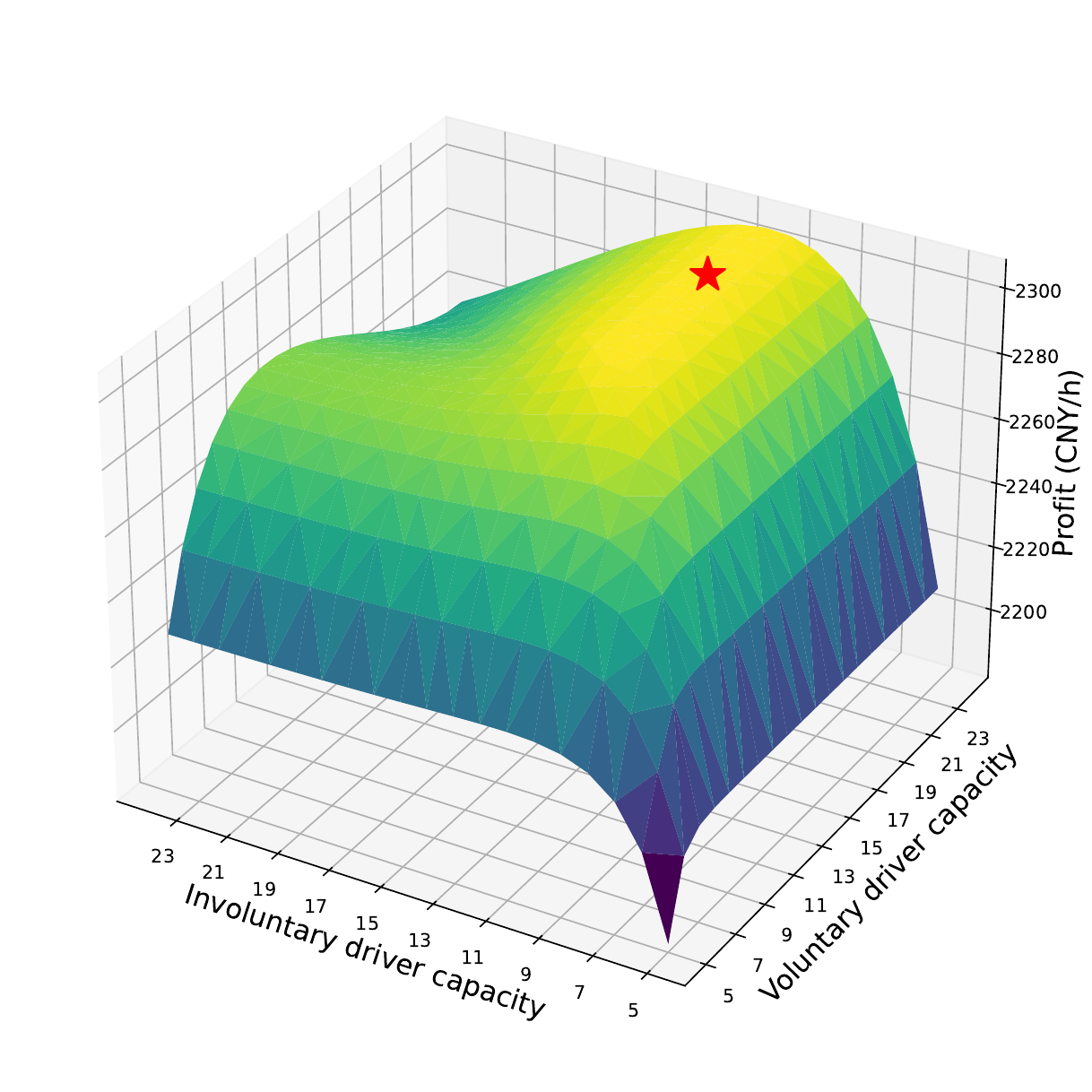}
  \caption{Surface plot for conditional maximum profit rate under different capacities}\label{fig:numer_PO_surface}
\end{figure}

We first solve the profit maximization problem (\ref{eq:profit_maximization_new}). For this case, the optimal capacity of the involuntary drivers and the voluntary drivers is 12 and 20, and the profit gain of the airport operations is 2,305.97 CNY per hour. The corresponding optimal lotteries for the involuntary and voluntary drivers are shown in  \Cref{fig:numer_MO_lotteries}. Then, the resultant equilibrium waiting times $\{W_m^n\}$ are shown in \Cref{fig:numer_PO_waittime}. Finally, for this example, we show the impact of different queue capacities on the lower-level objective function in \Cref{fig:numer_PO_surface}.

\Cref{fig:numer_MO_lotteries} demonstrates the optimal lotteries along with their properties. First, looking at the optimal lotteries, we can see that they always assign a higher probability to place an arriving driver at the front of the queue rather than the back of the queue. This is evident as the color of the grid turns warmer when the entry position gets smaller on the vertical axis, indicating that positions at the front are assigned a higher probability. Second, the tendency of placing the new arriving driver at the very first position of the queue decreases with the queue length. Looking at the bottom row, the color of the grid grows lighter when the queue length becomes larger, indicating that $\delta_m^{1,n}$ gets smaller when $n$ gets larger. Overall, the profit-maximizing lotteries is a policy between FIFO and last-in-first-out (LIFO). The shown lotteries are in one dimension similar to LIFO as they assign the highest probability to place new arrivals at the front of the queue, and on the other hand, it has some similarity to FIFO in the sense that the expected entry position of new arrivals increases with queue length.

\Cref{fig:numer_PO_waittime} reveals a key property of the profit-maximizing lottery: it ensures the equilibrium waiting time is constant, regardless of the queue length upon arrival. This phenomenon is a direct consequence of the platform's pricing strategy under its profit-maximization objective. The platform seeks to extract maximum revenue by setting a price at the limit of the drivers' participation constraint. This constraint is dictated by the disutility of waiting, meaning the maximum sustainable price is intrinsically linked to the maximum expected waiting time for each driver type. Therefore, by setting all expected waiting times equal, the platform can minimize the maximum expected waiting time and increase the price, therefore achieving profit maximization.

As illustrated in \Cref{fig:numer_PO_surface}, the profit function exhibits concavity with respect to the capacity variables, indicating that an optimal operating region exists. This optimum represents a critical trade-off: insufficient capacity limits service volume and thus suppresses revenue, whereas excessive capacity diminishes profitability by elevating system-wide waiting times and idle costs. The analysis further reveals that the objective function has a heightened sensitivity to the capacity allocated to involuntary drivers. This can be attributed to the higher arrival rate of this cohort; consequently, marginal changes to their capacity threshold induce more significant perturbations in the queueing system's steady-state probabilities and expected waiting times.

\begin{figure}[h!]
\centering
    \subfloat[][Involuntary drivers]{\includegraphics[width=0.5\linewidth]{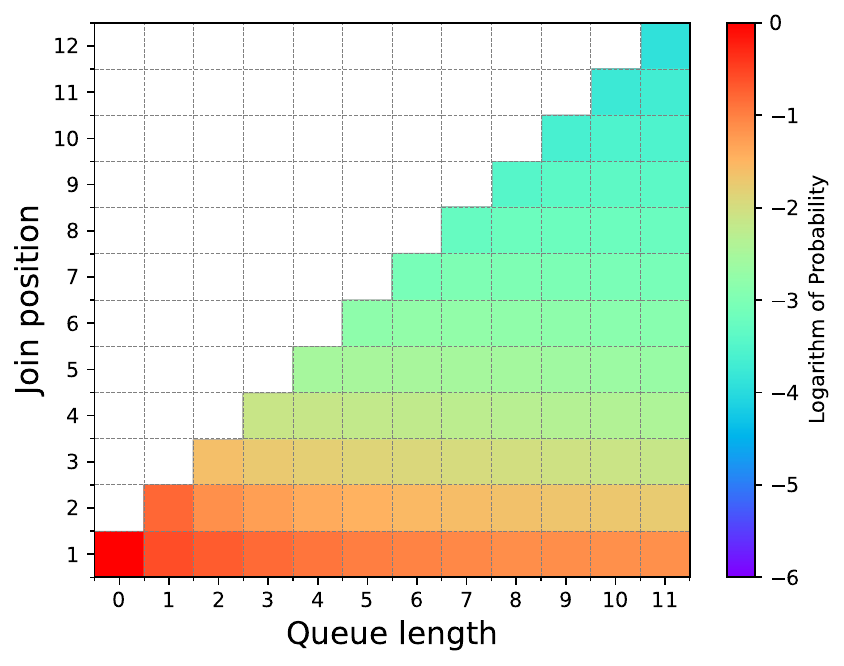}}
	\subfloat[][Voluntary drivers]{\includegraphics[width=0.5\linewidth]{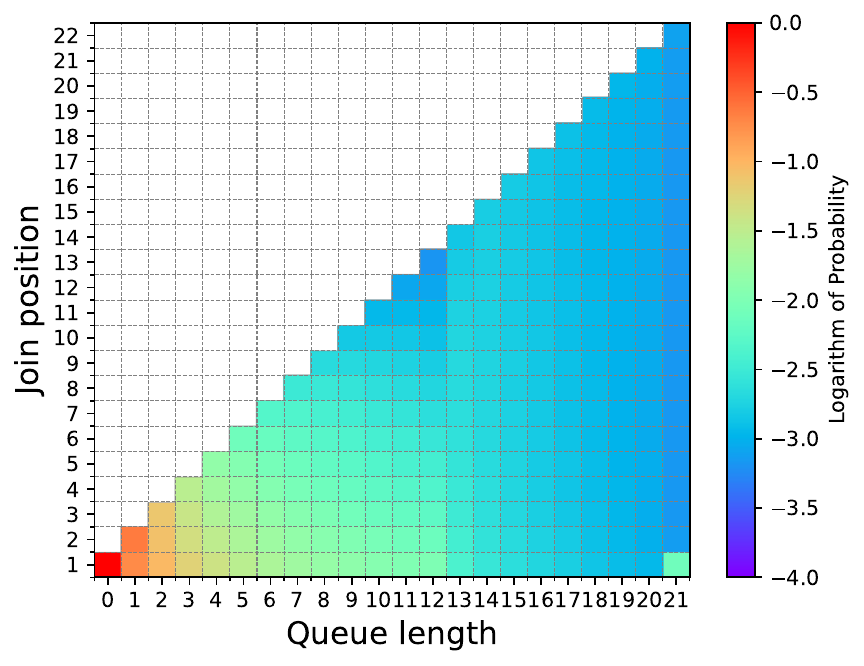}}
	\caption{The social-welfare-maximizing driver entry position lotteries}
	\label{fig:numer_SO_lotteries}
\end{figure}

\begin{figure}[h!]
  \centering
  \includegraphics[width=0.6\textwidth]{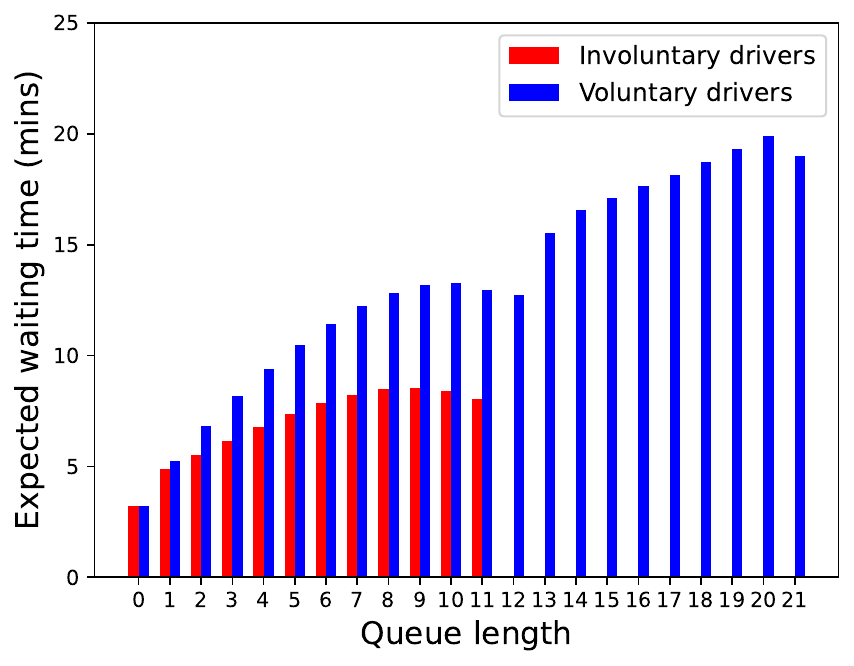}
  \caption{Equilibrium conditional waiting times for social welfare maximization}\label{fig:numer_SO_waittime}
\end{figure}

\begin{figure}[h!]
  \centering
  \includegraphics[width=0.6\textwidth]{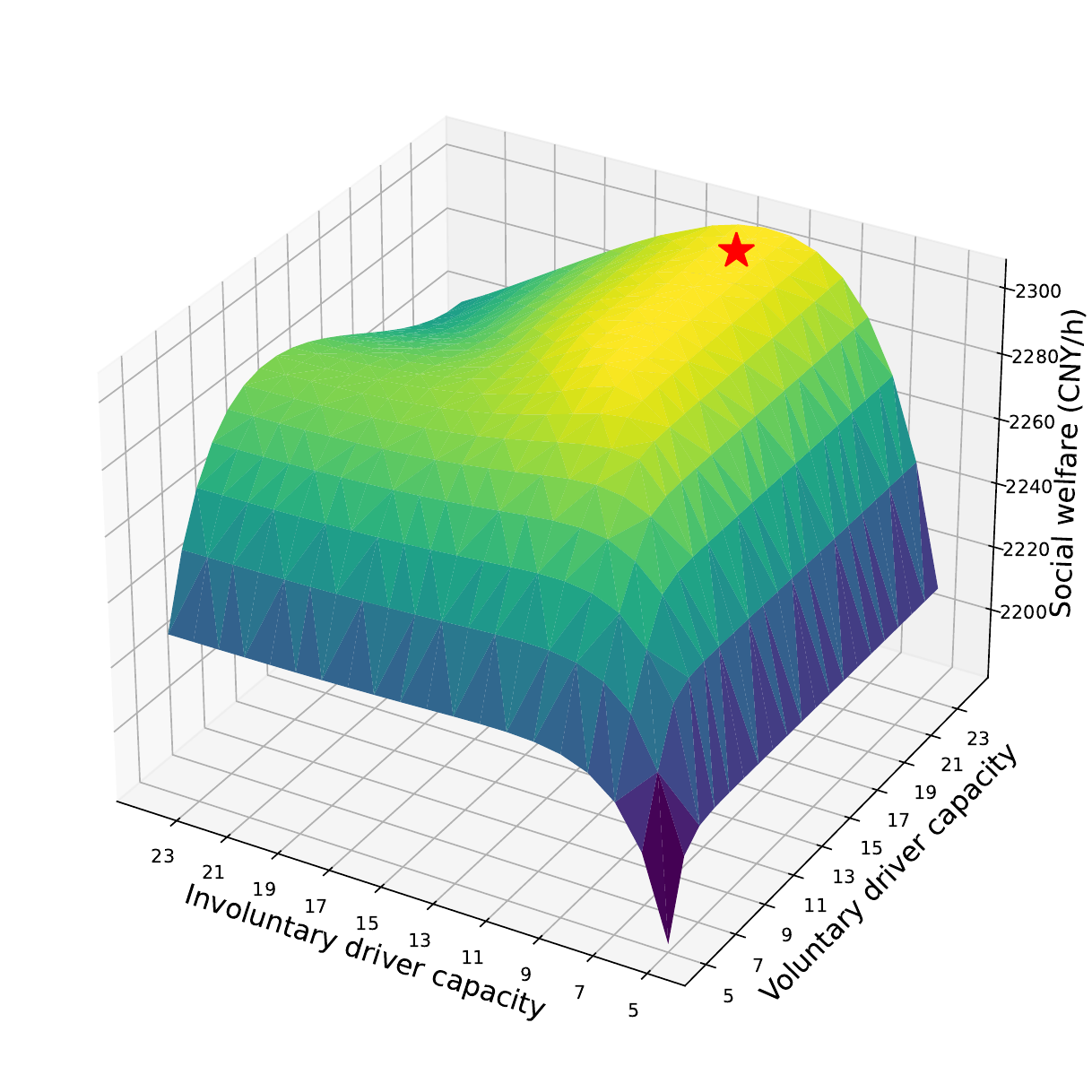}
  \caption{Surface plot for conditional maximum social welfare under different capacities}\label{fig:numer_SO_surface}
\end{figure}

For the social welfare maximization problem, under the same setting, the optimal capacities of the involuntary and voluntary drivers are 12 and 22. The social welfare gain of airport operations is 2,306 CNY per hour. The welfare-maximizing lotteries are shown in \Cref{fig:numer_SO_lotteries}, followed by the resultant equilibrium expected waiting times in \Cref{fig:numer_SO_waittime} and the impact of capacities on the maximum equilibrium social welfare in \Cref{fig:numer_SO_surface}.

The social-welfare-maximizing lotteries, depicted in \Cref{fig:numer_SO_lotteries}, exhibit structural properties analogous to the profit-maximizing lotteries shown in \Cref{fig:numer_MO_lotteries}. A primary characteristic is that for any given queue length, the probability mass of the optimal lottery is concentrated in front of the queue. However, this tendency diminishes as the queue length increases, with the probability of being placed at the front generally decreasing. Despite these similarities, there are significant distinctions driven by the different objective functions. First, the welfare-maximizing lottery for voluntary drivers (\Cref{fig:numer_SO_lotteries}(b)) displays a notable structural shift after the capacity of the other group, the involuntary drivers, has been reached. Beyond this point, the lottery distribute new arrivals more evenly across positions, suggesting a change in optimization strategy once a capacity constraint for one driver group becomes binding. Second, compared to their profit-maximizing counterparts, the social-welfare-maximizing lotteries give less priority to new arrivals, assigning a lower probability to being placed at the front of the queue and a correspondingly higher probability of being placed at the back. This is visually evident when comparing \Cref{fig:numer_MO_lotteries}(a) and \Cref{fig:numer_SO_lotteries}(a), where the latter shows a greater concentration of probability mass on the back of the queue.

In contrast to the profit-maximizing case, the expected waiting times under social welfare maximization are not constant, as shown in \Cref{fig:numer_SO_waittime}. While the expected waiting time generally increases with queue length, the plot reveals a non-monotonic behavior near the capacity threshold for involuntary drivers, where the waiting time dips slightly at queue lengths of 11 and 12. This decrement occurs because the system's dynamics change once this capacity limit is reached. When involuntary drivers no longer join, the overall rate of new arrivals --- any of whom could be placed ahead of existing drivers by the lottery --- is substantially reduced. This temporary reduction in the frequency of forward-queue insertions lowers the expected waiting time for those already in the queue, before the dominant trend of increasing wait times resumes.

\Cref{fig:numer_SO_surface} illustrates that properly calibrating capacities is crucial for maximizing social welfare; the system's performance declines when capacities are set too low or too high. Similarly, as shown in \Cref{fig:numer_PO_surface}, social welfare is more sensitive to the capacity of involuntary drivers than to that of voluntary drivers.

\subsection{Sensitivity analysis and control performance}

In this section, we provide additional numerical examples of lottery control under varying market conditions and regulations to illustrate: (1) how market conditions and regulations influence the structure and performance of optimal lottery control; and (2) the comparative performance of lottery control against dynamic and static pricing approaches.

\begin{figure}[h!]
\centering
    \subfloat[][Optimal capacities]{\includegraphics[width=0.5\linewidth]{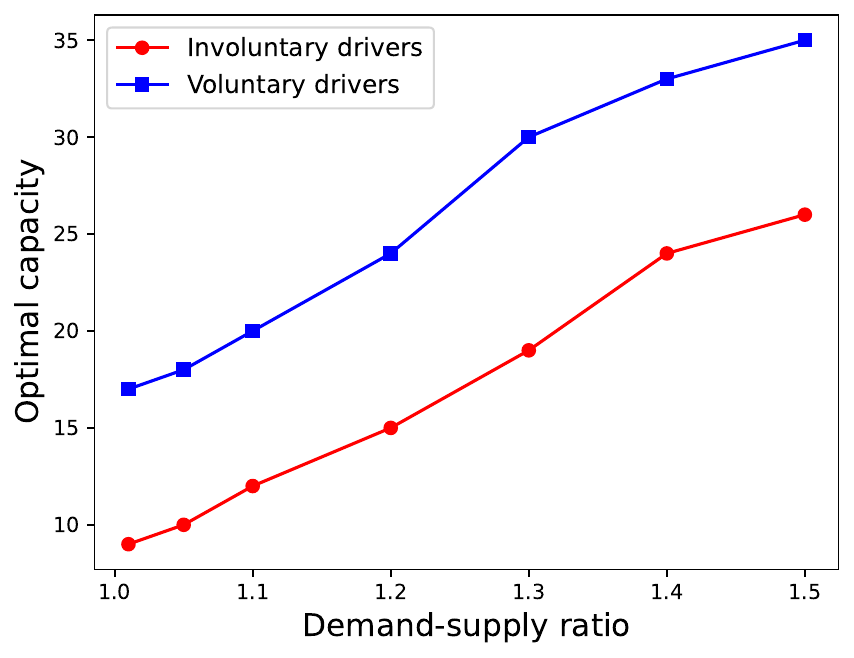}}
	\subfloat[][Average waiting time]{\includegraphics[width=0.5\linewidth]{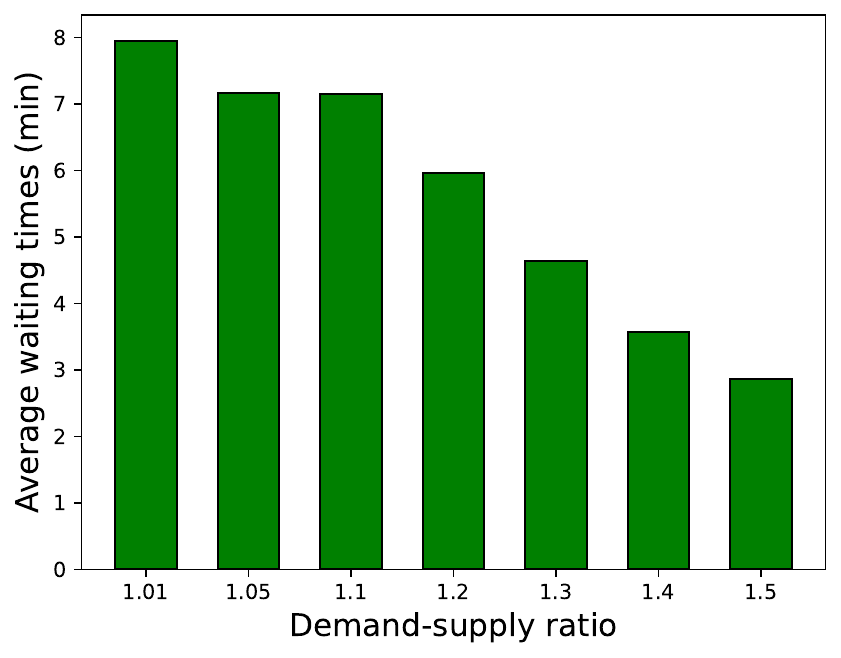}}
	\caption{Solution properties of profit-maximizing lottery control under different demand-supply ratios}
	\label{fig:numer_sensi_profit_capacities}
\end{figure}

\begin{figure}[h!]
  \centering
  \includegraphics[width=0.5\textwidth]{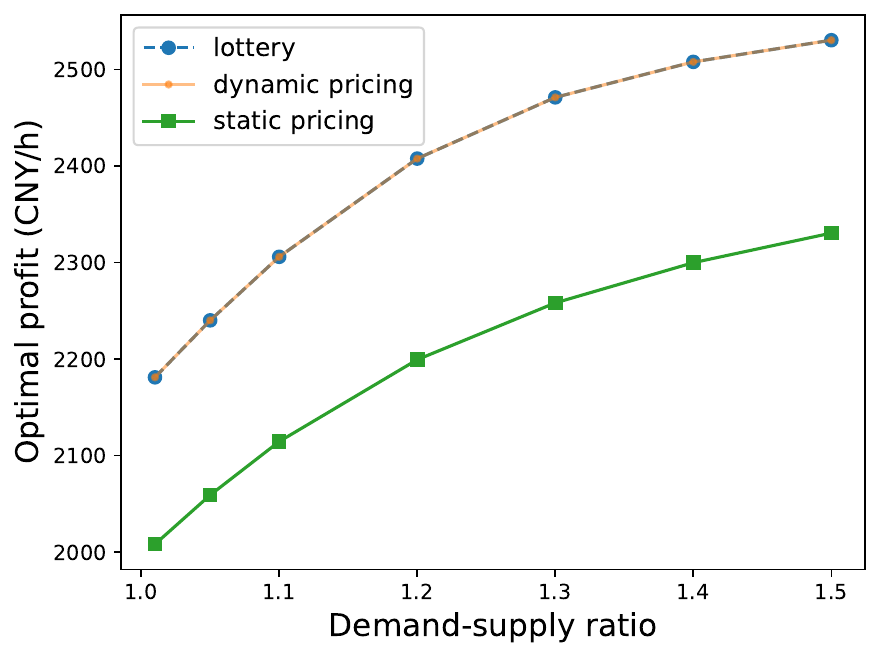}
  \caption{Optimal profits of the lottery control and pricing controls under different demand-supply ratios}\label{fig:numer_sensi_profit_performance}
\end{figure}

First, we analyze how market conditions affect the lottery control's structure and performance. We vary the overall demand-supply ratio, $\mu/(\lambda_1+\lambda_2)$, from 1.01 to 1.5 by adjusting the service rate $\mu$, with other parameters held as defined in \Cref{tb:para_model}. Under the profit-maximization objective, \Cref{fig:numer_sensi_profit_capacities} shows the resulting optimal capacities and average waiting times for both driver types. \Cref{fig:numer_sensi_profit_performance} then compares the profitability of our lottery control against two benchmarks: a FIFO policy with dynamic pricing and a FIFO policy with static pricing.

\Cref{fig:numer_sensi_profit_capacities} demonstrates that with the increment of demand, the passenger arrival rate increases and it is beneficial for the platform to allow more idle drivers into the virtual queue. This is shown in \Cref{fig:numer_sensi_profit_capacities}(a), as the optimal capacities for both voluntary drivers and involuntary drivers increases with the demand-supply ratio. Furthermore, as shown in \Cref{fig:numer_sensi_profit_capacities}(b), despite an increased number of drivers in the queue, the average waiting time of drivers decreases with the demand-supply ratio due to the more frequent arrival of passengers. Therefore, with the increase of passenger demand, the platform can maintain more drivers in the system while increasing the commission fee at the same time due to the reduced waiting time. This allows the platform to generate more profit from both the volume and pricing perspective.

As shown in \Cref{fig:numer_sensi_profit_performance}, the optimal profit for all three control methods increases with the demand-supply ratio, which is an intuitive result. However, the rate of this profit growth declines as the ratio increases, indicating diminishing marginal returns with respect to market tightness. The results also provide a clear performance comparison. Our proposed lottery control is highly effective, achieving profits nearly identical to the dynamic FIFO pricing benchmark. It significantly outperforms the static pricing policy, yielding an approximate 11\% increase in profit. Furthermore, this performance advantage over static pricing becomes more pronounced in high-demand scenarios, as the gap between the two methods widens with the increasing demand-supply ratio.

\begin{figure}[h!]
\centering
    \subfloat[][Optimal capacities]{\includegraphics[width=0.5\linewidth]{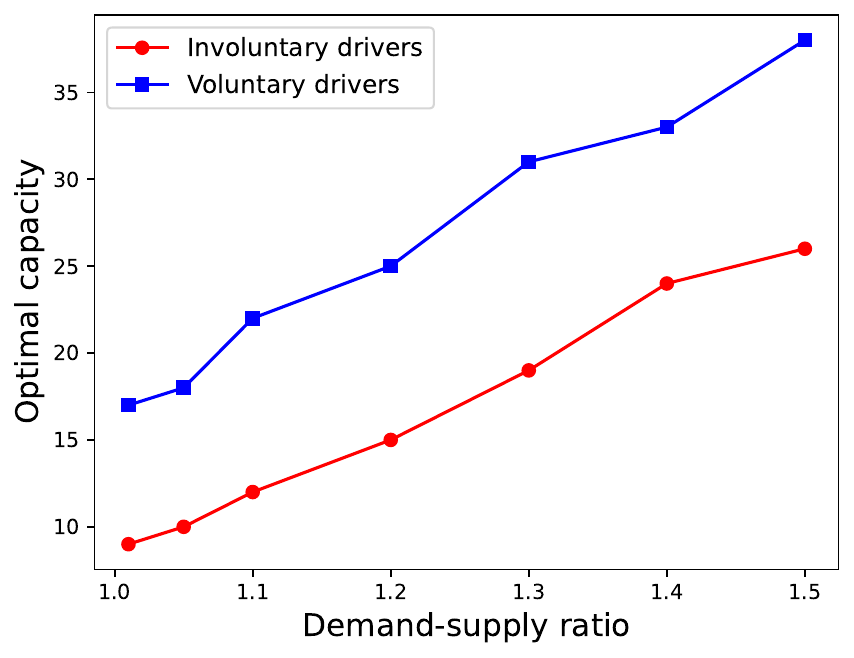}}
	\subfloat[][Average waiting time]{\includegraphics[width=0.5\linewidth]{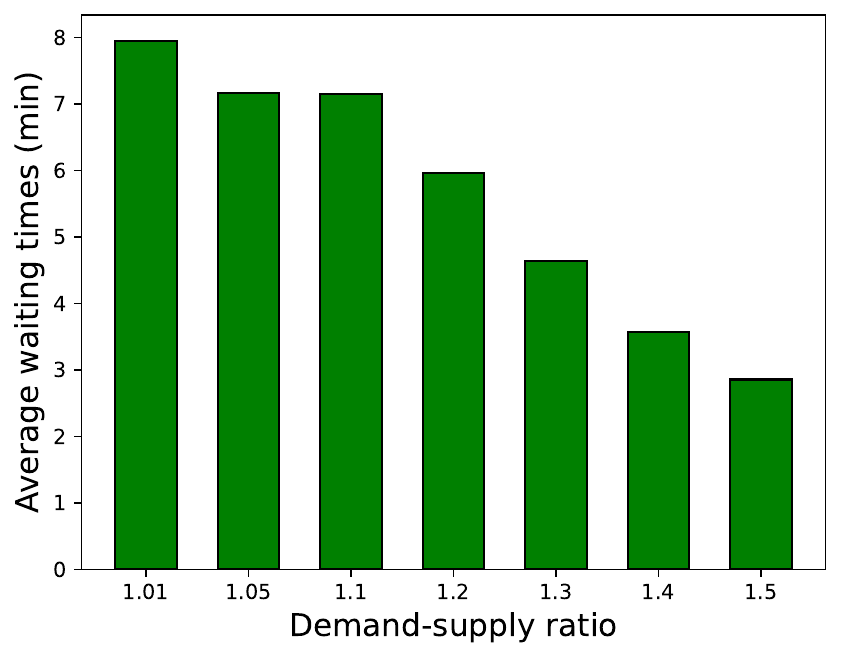}}
	\caption{Solution properties of social-welfare-maximizing lottery control under different demand-supply ratios}
	\label{fig:numer_sensi_welfare_capacities}
\end{figure}

\begin{figure}[h!]
  \centering
  \includegraphics[width=0.5\textwidth]{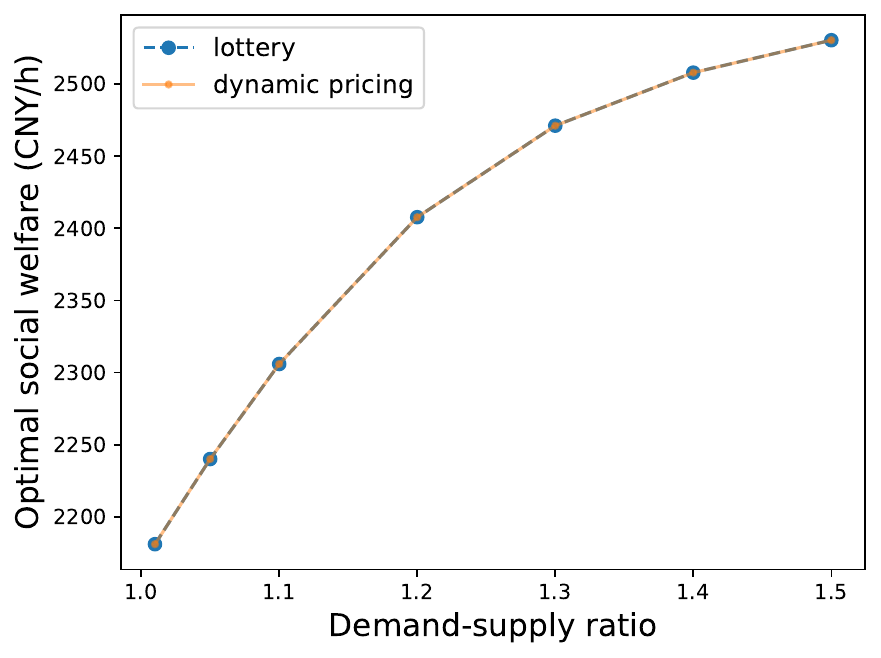}
  \caption{Optimal social welfare of the lottery control and dynamic pricing control under different demand-supply ratios}\label{fig:numer_sensi_welfare_performance}
\end{figure}

For the social welfare maximization problem, we present the solution properties of the lottery control in \Cref{fig:numer_sensi_welfare_capacities} and the relative performances of the lottery control and dynamic pricing in \Cref{fig:numer_sensi_welfare_performance}. We omit the performance of static pricing here as it can produce the same optimal social welfare as dynamic pricing \citep{hassin2016rational}.

The social-welfare-maximizing scenario, illustrated in \Cref{fig:numer_sensi_welfare_capacities}, exhibits trends analogous to the profit-maximizing case. As the demand-supply ratio increases, it becomes optimal to expand system capacity to accommodate more drivers. Notably, this capacity expansion is achieved without increasing the average waiting time, indicating an efficient scaling of the system. When compared to the profit-maximization outcomes, the optimal capacities under welfare maximization are largely consistent, though occasionally slightly higher. A more significant distinction is the reduction in average waiting times. Although this reduction is modest, its magnitude increases as the optimal capacity grows with the demand-supply ratio. Therefore, we anticipate that the welfare benefits of lottery control would be more pronounced in larger-scale queueing systems with greater operational flexibility.

As shown in \Cref{fig:numer_sensi_welfare_performance}, increasing the demand-supply ratio enhances social welfare, albeit with diminishing marginal returns—a trend consistent with the profit-maximization case. In terms of performance, lottery control and dynamic pricing yield nearly identical social welfare outcomes. When combined with previous findings (\Cref{fig:numer_sensi_profit_performance}), this result demonstrates the robustness of lottery control, as it consistently matches the performance of the dynamic pricing benchmark across both profit and welfare objectives.

\begin{figure}[h!]
  \centering
  \includegraphics[width=0.5\textwidth]{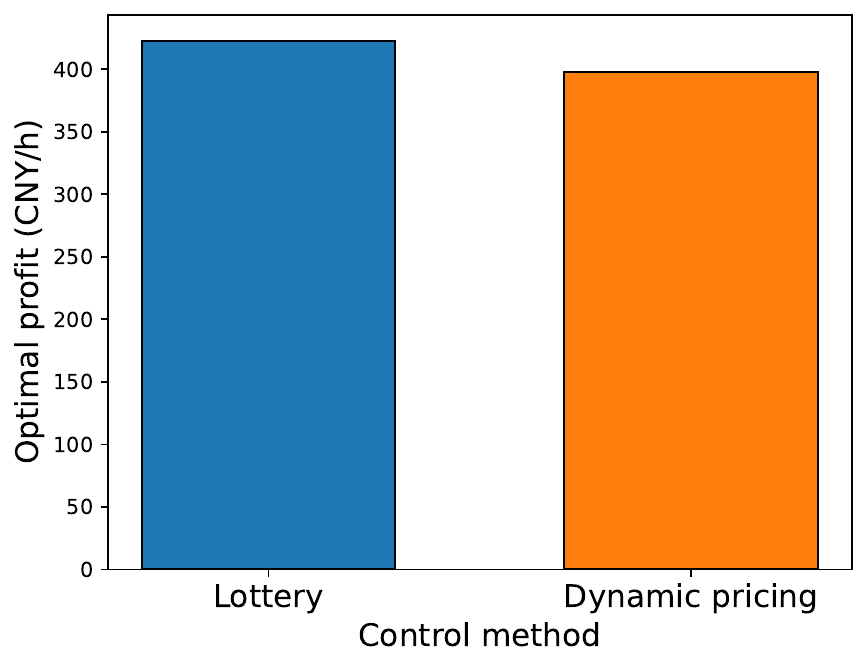}
  \caption{Comparision of lottery control and dynamic pricing control under a commission cap}\label{fig:numer_sensi_cap_compare}
\end{figure}

Commission fee caps are a popular regulatory tool in ride-sourcing markets \citep{zha2016economic,vignon2021regulating}. In such a regulated environment, we demonstrate that lottery control can outperform the dynamic pricing benchmark. We illustrate this advantage by modeling a two-group system under the commission fee cap policy defined in \Cref{tb:para_model_cap}. The resulting comparison of optimal profits under both control schemes is presented in \Cref{fig:numer_sensi_cap_compare}.

\begin{table}[h!]\centering
\caption{Value of commission fee cap market model parameters}
\label{tb:para_model_cap}\footnotesize
\begin{tabular}{c c}
\hline
Parameter & Value\\
\hline
$\lambda_1$ & 93.8 veh/h\\
$\lambda_2$  & 31.9 veh/h\\
$\mu$ & 132.0 orders/h\\
$R_1$  & 7.5 CNY\\
$R_2$  & 12.5 CNY\\
$r$  & 40 CNY/h\\
$v$  & 10 CNY/h\\
$T_d$ & 15 minutes\\
$p_{max}$ & 4.25 CNY\\
\hline
\end{tabular}
\end{table}

\Cref{fig:numer_sensi_cap_compare} confirms that lottery control can achieve superior profitability over dynamic pricing when the market is regulated by a commission fee cap. This advantage arises because the two methods have different control levers. A price cap directly constrains the primary tool of a dynamic pricing system, limiting its ability to optimize. In contrast, lottery control leverages queue position to manage driver utility and wait times, allowing it to find a more globally optimal solution under the regulatory constraint.


\section{Conclusion} \label{sec:conclusion}
In this paper, we study the idle ride-sourcing driver management problem at transportation terminals in a virtual queue setting. We model driver behavior and system dynamics through a queueing system and introduce a novel lottery-based approach that implicitly guides idle drivers' balking decisions by influencing their queue-length-conditional expected waiting times through controlled queue entry positions.

Using this model, we formulate both profit and social welfare maximization problems to determine the optimal lottery control. These problems are initially intractable, manifesting as infinite-dimensional optimization with non-smooth equilibrium constraints. To overcome this, we leverage a key theoretical insight: the optimal equilibrium exhibits a threshold structure, where drivers join below a certain capacity and balk above it. This property allows us to reformulate the problem into a tractable bi-level program by decomposing the equilibrium, which resolves the non-smoothness. On such transformation's properties, we further prove that this reformulation preserves the optimality of the social welfare solution.

Our analysis yields several key findings. Theoretically, we prove that lottery control outperforms dynamic pricing in maximizing social welfare, due to its fine-grained ability to manage individual waiting times and minimize system-wide idle costs. Numerically, we find that for profit maximization, our method's performance is comparable to dynamic pricing and improves upon static pricing by approximately 11\%. We also demonstrate that lottery control can surpass dynamic pricing under regulatory constraints, such as a commission fee cap, highlighting its operational robustness. These combined findings establish lottery control as a competitive alternative to dynamic pricing levers for supply management.

While our approach provides valuable insights into idle ride-sourcing driver management, several limitations remain for future research. First, although our numerical analysis indicates that the optimality gap between lottery control and dynamic pricing is small, a theoretical characterization of this gap remains unknown. Future work can explore the formal properties of the optimality gap and develop new heuristic approaches that relax the structural restrictions imposed by \Cref{prop:threshold_all_join_profit}. Second, our analysis primarily focuses on driver-side management by adjusting commission fees and lotteries, without explicitly modeling passenger behavior. Future research can incorporate passenger-side dynamics, allowing for the analysis of lottery control across both market sides and its impact on system performance. Finally, we adopt the classical rational queueing model from \cite{palm1957research}, which does not account for reneging behavior, as its inclusion introduces significant complications on steady-state analysis. Future research can explore methods to incorporate reneging behavior into the system and examine the resulting optimal lottery structures and their performance implications.
\bibliographystyle{apalike} 

\begin{thebibliography}{}

\bibitem[Adarkwah, 2023]{adarkwah2023reducing}
Adarkwah, B.~K. (2023).
\newblock Reducing excessive dwell time of airport queueing lot users at o’hare international airport.
\newblock Master's thesis, University of Illinois at Chicago.

\bibitem[Agarwal et~al., 2025]{agarwal2023information}
Agarwal, S., Cheng, S.-F., Keppo, J., Wang, L., and Yang, Z. (2025).
\newblock Information provision and labor efficiency: Evidence from taxis in singapore.
\newblock {\em Available at SSRN 4234268}.

\bibitem[Alnaggar et~al., 2024]{alnaggar2024heatmap}
Alnaggar, A., Gzara, F., and Bookbinder, J. (2024).
\newblock Heatmap design for probabilistic driver repositioning in crowdsourced delivery.
\newblock {\em Transportation Science}.

\bibitem[Ashkrof et~al., 2020]{ashkrof2020understanding}
Ashkrof, P., de~Almeida~Correia, G.~H., Cats, O., and Van~Arem, B. (2020).
\newblock Understanding ride-sourcing drivers' behaviour and preferences: Insights from focus groups analysis.
\newblock {\em Research in Transportation Business \& Management}, 37:100516.

\bibitem[{Baltimore Regional Transportation Board}, 2023]{BALreport}
{Baltimore Regional Transportation Board} (2023).
\newblock 2023 {W}ashington-{B}altimore regional air passenger survey.

\bibitem[Beojone and Geroliminis, 2023]{beojone2023relocation}
Beojone, C.~V. and Geroliminis, N. (2023).
\newblock Relocation incentives for ride-sourcing drivers with path-oriented revenue forecasting based on a markov chain model.
\newblock {\em Transportation Research Part C: Emerging Technologies}, 157:104375.

\bibitem[Besbes et~al., 2021]{besbes2021surge}
Besbes, O., Castro, F., and Lobel, I. (2021).
\newblock Surge pricing and its spatial supply response.
\newblock {\em Management Science}, 67(3):1350--1367.

\bibitem[Bimpikis et~al., 2019]{bimpikis2019spatial}
Bimpikis, K., Candogan, O., and Saban, D. (2019).
\newblock Spatial pricing in ride-sharing networks.
\newblock {\em Operations Research}, 67(3):744--769.

\bibitem[Chen et~al., 2021]{chen2021spatial}
Chen, C., Yao, F., Mo, D., Zhu, J., and Chen, X.~M. (2021).
\newblock Spatial-temporal pricing for ride-sourcing platform with reinforcement learning.
\newblock {\em Transportation Research Part C: Emerging Technologies}, 130:103272.

\bibitem[Chen and Frank, 2001]{chen2001state}
Chen, H. and Frank, M.~Z. (2001).
\newblock State dependent pricing with a queue.
\newblock {\em IIE Transactions}, 33(10):847--860.

\bibitem[Chen et~al., 2019]{chen2019inbede}
Chen, H., Jiao, Y., Qin, Z., Tang, X., Li, H., An, B., Zhu, H., and Ye, J. (2019).
\newblock Inbede: Integrating contextual bandit with td learning for joint pricing and dispatch of ride-hailing platforms.
\newblock In {\em 2019 IEEE International Conference on Data Mining (ICDM)}, pages 61--70. IEEE.

\bibitem[Freund and van Ryzin, 2021]{freund2021pricing}
Freund, D. and van Ryzin, G. (2021).
\newblock Pricing fast and slow: Limitations of dynamic pricing mechanisms in ride-hailing.
\newblock {\em Available at SSRN 3931844}.

\bibitem[Guda and Subramanian, 2019]{guda2019your}
Guda, H. and Subramanian, U. (2019).
\newblock Your uber is arriving: Managing on-demand workers through surge pricing, forecast communication, and worker incentives.
\newblock {\em Management Science}, 65(5):1995--2014.

\bibitem[Haferkamp et~al., 2024]{haferkamp2024heatmap}
Haferkamp, J., Ulmer, M.~W., and Ehmke, J.~F. (2024).
\newblock Heatmap-based decision support for repositioning in ride-sharing systems.
\newblock {\em Transportation Science}, 58(1):110--130.

\bibitem[Hassin, 2016]{hassin2016rational}
Hassin, R. (2016).
\newblock {\em Rational queueing}.
\newblock CRC press.

\bibitem[Hermawan and Regan, 2018]{hermawan2018impacts}
Hermawan, K. and Regan, A.~C. (2018).
\newblock Impacts on vehicle occupancy and airport curb congestion of transportation network companies at airports.
\newblock {\em Transportation Research Record}, 2672(23):52--58.

\bibitem[Hu et~al., 2022]{hu2022surge}
Hu, B., Hu, M., and Zhu, H. (2022).
\newblock Surge pricing and two-sided temporal responses in ride hailing.
\newblock {\em Manufacturing \& Service Operations Management}, 24(1):91--109.

\bibitem[Ji and Cheng, 2021]{ji2021automated}
Ji, M. and Cheng, S.-F. (2021).
\newblock Automated taxi queue management at high-demand venues.
\newblock In {\em 2021 IEEE 17th International Conference on Automation Science and Engineering (CASE)}, pages 1757--1762. IEEE.

\bibitem[Leiner and Adler, 2020]{leiner2020transportation}
Leiner, C. and Adler, T. (2020).
\newblock {\em Transportation network companies (TNCs): impacts to airport revenues and operations—reference guide}.
\newblock Number Project 01-35.

\bibitem[{L.E.K Consulting}, 2018]{SFOreport}
{L.E.K Consulting} (2018).
\newblock The future of airport ground access.

\bibitem[Liu et~al., 2025]{liu2025bounded}
Liu, T., Xu, Z., Keepo, J., Yin, Y., and Zhu, H. (2025).
\newblock Bounded rationality in ride-sourcing drivers' dwelling at transportation terminals: A behavioral queueing analysis.
\newblock {\em Available at SSRN 5050040}.

\bibitem[{Lyft}, 2024]{Lyftairport}
{Lyft} (2024).
\newblock Receiving airport fifo pickup requests.
\newblock url = {https://help.lyft.com/hc/en-us/all/articles/115012922787-Receiving-Airport-FIFO-pickup-requests}.

\bibitem[{Mission Local}, 2023]{SFOqueue}
{Mission Local} (2023).
\newblock Uber, lyft drivers at sfo doing lots of waiting, less driving.
\newblock url = {https://missionlocal.org/2023/03/uber-lyft-drivers-sfo-lots-of-waiting/}.

\bibitem[{New York City Taxi and Limousine Commission}, 2022]{NYCHearing}
{New York City Taxi and Limousine Commission} (2022).
\newblock Notice of public hearing and opportunity to comment on proposed rules.
\newblock url = {https://www.nyc.gov/assets/tlc/downloads/pdf/proposed-rules-04-5-2022.pdf}.

\bibitem[Nourinejad and Ramezani, 2020]{nourinejad2020ride}
Nourinejad, M. and Ramezani, M. (2020).
\newblock Ride-sourcing modeling and pricing in non-equilibrium two-sided markets.
\newblock {\em Transportation Research Part B: Methodological}, 132:340--357.

\bibitem[{\"O}zkan, 2020]{ozkan2020joint}
{\"O}zkan, E. (2020).
\newblock Joint pricing and matching in ride-sharing systems.
\newblock {\em European Journal of Operational Research}, 287(3):1149--1160.

\bibitem[Palm, 1957]{palm1957research}
Palm, R. C.~A. (1957).
\newblock {\em Research on telephone traffic carried by full availability groups}.
\newblock Tele.

\bibitem[Pollio, 2021]{pollio2021uber}
Pollio, A. (2021).
\newblock Uber, airports, and labour at the infrastructural interfaces of platform urbanism.
\newblock {\em Geoforum}, 118:47--55.

\bibitem[Sadeghi and Smith, 2019]{sadeghi2019re}
Sadeghi, A. and Smith, S.~L. (2019).
\newblock On re-balancing self-interested agents in ride-sourcing transportation networks.
\newblock In {\em 2019 IEEE 58th Conference on Decision and Control (CDC)}, pages 5119--5125. IEEE.

\bibitem[Seele et~al., 2021]{seele2021mapping}
Seele, P., Dierksmeier, C., Hofstetter, R., and Schultz, M.~D. (2021).
\newblock Mapping the ethicality of algorithmic pricing: A review of dynamic and personalized pricing.
\newblock {\em Journal of Business Ethics}, 170:697--719.

\bibitem[{Uber}, 2024]{UberHKG}
{Uber} (2024).
\newblock Driving at {Hong} {Kong} international airport ({HKG}).
\newblock url = https://www.uber.com/global/en/r/airports/hkg/pickup/.

\bibitem[Varma et~al., 2023]{varma2023dynamic}
Varma, S.~M., Bumpensanti, P., Maguluri, S.~T., and Wang, H. (2023).
\newblock Dynamic pricing and matching for two-sided queues.
\newblock {\em Operations Research}, 71(1):83--100.

\bibitem[Vignon et~al., 2021]{vignon2021regulating}
Vignon, D.~A., Yin, Y., and Ke, J. (2021).
\newblock Regulating ridesourcing services with product differentiation and congestion externality.
\newblock {\em Transportation Research Part C: Emerging Technologies}, 127:103088.

\bibitem[Wadud, 2020]{wadud2020examination}
Wadud, Z. (2020).
\newblock An examination of the effects of ride-hailing services on airport parking demand.
\newblock {\em Journal of Air Transport Management}, 84:101783.

\bibitem[Wang et~al., 2025]{wangexpress}
Wang, Y., Lin, X., He, F., Xu, Z., and Shen, Z.-J.~M. (2025).
\newblock Express: Efficiency and interventions of strategic driver relocation for ride-hailing platforms.
\newblock {\em Production and Operations Management}, page 10591478251376798.

\bibitem[Xu, 2020]{xu2020empty}
Xu, Z. (2020).
\newblock {\em On the Empty Miles of Ride-Sourcing Services: Theory, Observation and Countermeasures}.
\newblock PhD thesis, University of Michigan.

\bibitem[Xu et~al., 2020]{xu2020recommender}
Xu, Z., Men, C., Li, P., Jin, B., Li, G., Yang, Y., Liu, C., Wang, B., and Qie, X. (2020).
\newblock When recommender systems meet fleet management: Practical study in online driver repositioning system.
\newblock In {\em Proceedings of the web conference 2020}, pages 2220--2229.

\bibitem[Xu et~al., 2021]{xu2021generalized}
Xu, Z., Yin, Y., Chao, X., Zhu, H., and Ye, J. (2021).
\newblock A generalized fluid model of ride-hailing systems.
\newblock {\em Transportation Research Part B: Methodological}, 150:587--605.

\bibitem[Yang and Ying, 2024]{yang2024learning}
Yang, Z. and Ying, L. (2024).
\newblock Learning-based pricing and matching for two-sided queues.
\newblock {\em arXiv preprint arXiv:2403.11093}.

\bibitem[Yengejeh and Smith, 2021]{yengejeh2021rebalancing}
Yengejeh, A.~S. and Smith, S.~L. (2021).
\newblock Rebalancing self-interested drivers in ride-sharing networks to improve customer wait-time.
\newblock {\em IEEE Transactions on Control of Network Systems}, 8(4):1637--1648.

\bibitem[Yu et~al., 2022]{yu2022delay}
Yu, Q., Zhang, Y., and Zhou, Y.-P. (2022).
\newblock Delay information in virtual queues: A large-scale field experiment on a major ride-sharing platform.
\newblock {\em Management Science}, 68(8):5745--5757.

\bibitem[Zha et~al., 2018]{zha2018surge}
Zha, L., Yin, Y., and Du, Y. (2018).
\newblock Surge pricing and labor supply in the ride-sourcing market.
\newblock {\em Transportation Research Part B: Methodological}, 117:708--722.

\bibitem[Zha et~al., 2016]{zha2016economic}
Zha, L., Yin, Y., and Yang, H. (2016).
\newblock Economic analysis of ride-sourcing markets.
\newblock {\em Transportation Research Part C: Emerging Technologies}, 71:249--266.

\bibitem[Zhang et~al., 2024]{zhang2024dynamic}
Zhang, J., Hu, L., Li, Y., Xu, W., and Jiang, Y. (2024).
\newblock Dynamic joint decision of matching parameters and relocation strategies in ride-sourcing systems interacting with traffic congestion.
\newblock {\em Transportation Research Part C: Emerging Technologies}, 161:104524.

\bibitem[Zhu et~al., 2021]{zhu2021mean}
Zhu, Z., Ke, J., and Wang, H. (2021).
\newblock A mean-field markov decision process model for spatial-temporal subsidies in ride-sourcing markets.
\newblock {\em Transportation Research Part B: Methodological}, 150:540--565.

\end{thebibliography}


\end{document}